\documentclass[oneside, reqno] {amsart}

\usepackage{xypic}
\input xy
\xyoption{all}
\usepackage{epsfig}
\usepackage{amsthm}
\usepackage{amssymb}
\usepackage{amsmath}
\usepackage{amscd}
\usepackage{color}
\usepackage[T1]{fontenc}
\usepackage{upgreek}
\usepackage[font=scriptsize]{caption}
\usepackage{wrapfig}
\usepackage{graphicx}
\usepackage{subfigure}

\newcommand{\bg}{\begin{equation}}
\newcommand{\ed}{\end{equation}}
\newcommand{\bga}{\begin{eqnarray}}
\newcommand{\eda}{\end{eqnarray}}
\newcommand{\pf}{\textbf{Proof:\ }}

\def\cbdu{\par{\raggedleft$\Box$\par}}

\newtheorem {Theorem}  {Theorem}

\numberwithin{Theorem}{section}

\newtheorem {Lemma}[Theorem]  {Lemma}
\newtheorem {Proposition}[Theorem]{Proposition}
\theoremstyle{definition}
\newtheorem{Definition}[Theorem]{Definition}
\theoremstyle{remark}

\expandafter\chardef\csname pre amssym.def
at\endcsname=\the\catcode`\@ \catcode`\@=11
\def\undefine#1{\let#1\undefined}
\def\newsymbol#1#2#3#4#5{\let\next@\relax
 \ifnum#2=\@ne\let\next@\msafam@\else
 \ifnum#2=\tw@\let\next@\msbfam@\fi\fi
 \mathchardef#1="#3\next@#4#5}
\def\mathhexbox@#1#2#3{\relax
 \ifmmode\mathpalette{}{\m@th\mathchar"#1#2#3}%
 \else\leavevmode\hbox{$\m@th\mathchar"#1#2#3$}\fi}
\def\hexnumber@#1{\ifcase#1 0\or 1\or 2\or 3\or 4\or 5\or 6\or 7\or 8\or
 9\or A\or B\or C\or D\or E\or F\fi}

\font\teneufm=eufm10 \font\seveneufm=eufm7 \font\fiveeufm=eufm5
\newfam\eufmfam
\textfont\eufmfam=\teneufm \scriptfont\eufmfam=\seveneufm
\scriptscriptfont\eufmfam=\fiveeufm

\catcode`\@=\csname pre amssym.def at\endcsname

\newcounter{remark}
\def  \12  {{\frac{1}{2}}}

\def\build#1_#2^#3{\mathrel{\mathop{\kern 0pt#1}\limits_{#2}^{#3}}}

\numberwithin{equation}{section}

\begin{document}
%\currannalsline{0}{2006}

%\title[Local existence for Hall-MHD]{Local existence for the non-resistive Hall-magneto-hydrodynamics system in $\R^n$}
\title[1D model for MHD]{1D model for the 3D magnetohydrodynamics}

%\author{hello}

\author [Mimi Dai]{Mimi Dai}

\address{Department of Mathematics, Statistics and Computer Science, University of Illinois at Chicago, Chicago, IL 60607, USA}
\email{mdai@uic.edu} 

\author [Bhakti Vyas]{Bhakti Vyas}

\address{Department of Mathematics, Statistics and Computer Science, University of Illinois at Chicago, Chicago, IL 60607, USA}
\email{bvyas2@uic.edu}

\author [Xiangxiong Zhang]{Xiangxiong Zhang}

\address{Department of Mathematics, Purdue University, West Lafayette, IN 47907, USA}
\email{zhan1966@purdue.edu}

\thanks{M. Dai and B. Vyas are partially supported by the NSF grants DMS--1815069 and DMS--2009422. X. Zhang is partially supported by the NSF grant DMS--1913120.}

\begin{abstract}

We propose a one-dimensional (1D) model for the three-dimensional (3D) incompressible ideal magnetohydrodynamics.  We establish a regularity criterion of the Beale-Kato-Majda type for this 1D model. Without the stretching effect, the model with only transport effect equipped with a proper sign is shown to have global in time strong solution. Some numerical simulations suggest that solutions of this model with smooth periodic initial data do not tend to develop singularities at finite time.

\bigskip


\hspace{0.02cm}CLASSIFICATION CODE: 35Q35, 65M70, 76D03, 76W05, 76B03.
\end{abstract}

\maketitle

\section{Introduction}
The ideal incompressible magnetohydrodynamics (MHD) governed by the set of partial differential equations% with Hall effect
\begin{equation}\label{mhd}
\begin{split}
u_t+(u\cdot\nabla) u-(B\cdot\nabla) B+\nabla \Pi=&\ 0, \\
B_t+(u\cdot\nabla) B-(B\cdot\nabla) u =&\ 0, \\
%+d_i \nabla\times ((\nabla\times B)\times B)=&\ 0, \\
\nabla\cdot u=0, \ \ \nabla\cdot B=&\ 0,
\end{split}
\end{equation}
is an important model in geophysics and astrophysics. In the system, the vector fields $u$ and $B$ denote the fluid velocity and magnetic field respectively; the scalar function $\Pi$ is the pressure. 
 We notice that (\ref{mhd}) reduces to the incompressible Euler equation if $B\equiv0$,
\begin{equation}\label{euler}
\begin{split}
u_t+(u\cdot\nabla) u+\nabla \Pi=&\ 0, \\
\nabla\cdot u=&\ 0.
\end{split}
\end{equation}
The mathematical question of whether or not a solution of the 3D Euler (\ref{euler}) develops singularity at finite time remains open. So does it for the 3D MHD (\ref{mhd}). 

Denote the vorticity by $\omega=\nabla\times u$. Taking a curl on (\ref{euler}) gives 
\begin{subequations}
\begin{align}
\omega_t+(u\cdot\nabla) \omega+(\omega\cdot\nabla) u=&\ 0, \label{vorticity}\\
u=&\ \nabla\times(-\Delta)^{-1}\omega.  \label{BS-eu}
\end{align}
\end{subequations}
We note that $u$ can be recovered from $\omega$ through the Biot-Savart law (\ref{BS-eu}) which involves a nonlocal operator. In (\ref{vorticity}), the quadratic term $(u\cdot\nabla) \omega$ is regarded as the transport term, while $(\omega\cdot\nabla) u$ represents the stretching effect. The general belief is that the stretching effect is responsible for dramatic wild behaviours of solutions, for instance, the appearance of finite-time singularity. 

\medskip

\subsection{1D models for Euler equation and related equations}
To gain insights towards understanding the properties of solutions to the Euler equation (\ref{euler}), approximating models and toy models have been proposed and studied in the literature. One type of 1D models for the vorticity form of Euler equation has attracted a great deal of attention, which can be traced back to the work of Constantin, Lax and Majda \cite{CLM}. The authors of \cite{CLM} proposed the following 1D model for system (\ref{vorticity})-(\ref{BS-eu}),
\begin{subequations}
\begin{align}
\omega_t=&\ \omega H\omega, \label{eq1-clm}\\
u_x=&\ H\omega, \label{eq2-clm}
\end{align}
\end{subequations}
with $\omega=\omega(t,x)$ and $u=u(t,x)$ for $t\geq 0$ and $x\in\mathbb R$. In the system, $H$ denotes the Hilbert transform defined by
\begin{equation}\label{H1}
H f=\frac{1}{\pi} P.V. \int_{-\infty}^{\infty}\frac{f(y)}{x-y}\, dy.
\end{equation}
We note that equation (\ref{eq2-clm}) is a 1D analogue of the Biot-Savart law (\ref{BS-eu}). With only stretching effect in equation (\ref{eq1-clm}), the authors solved system (\ref{eq1-clm})-(\ref{eq2-clm}) exactly and showed the formation of finite-time singularities for a class of initial data. Since then, various generalisations of (\ref{eq1-clm})-(\ref{eq2-clm}) have been studied both analytically and numerically. The De Gregorio model \cite{DeG1, DeG2}
\begin{subequations}
\begin{align}
\omega_t+u\omega_x- \omega H\omega=&\ 0, \label{eq1-deg}\\
u_x=&\ H\omega, \label{eq2-deg}
\end{align}
\end{subequations}
includes both transport and stretching effects. Numerical results of \cite{DeG1, DeG2} provide evidence that finite-time blow-up may not occur for system (\ref{eq1-deg})-(\ref{eq2-deg}). It indicates that the convection (transport) term has a regularization effect and it dominates the stretching term. Later on, in order to understand the competing effects of convection and stretching terms, Okamoto, Sakajo and Wunsch \cite{OSW} suggested to study the following family of models
\begin{subequations}
\begin{align}
\omega_t+au\omega_x- \omega H\omega=&\ 0, \label{eq1-osw}\\
u_x=&\ H\omega, \label{eq2-osw}
\end{align}
\end{subequations}
with a parameter $a\in \mathbb R$. The authors also conjectured global in time existence of solutions to (\ref{eq1-osw})-(\ref{eq2-osw}) with $a=1$ which is the De Gregorio model (\ref{eq1-deg})-(\ref{eq2-deg}). Indeed, Jia, Stewart and \v{S}ver\'ak \cite{Jia2019} proved that solutions of (\ref{eq1-deg})-(\ref{eq2-deg}) with initial data near a steady state are global and converge to this steady state. In contrast, Elgindi and Jeong \cite{EJ} showed singularity formation for (\ref{eq1-deg})-(\ref{eq2-deg}) in classes of H\"older continuous solutions. Moreover, the authors of \cite{EJ} established that, there exists smooth initial data such that solution of the Okamoto-Sakajo-Wunsch model (\ref{eq1-osw})-(\ref{eq2-osw}) with small $|a|$ develops self-similar type of blow-up at finite time. Later on, Elgindi, Ghoul and Masmoudi \cite{EGM} further showed that such self-similar blow-up is stable. 

When $a=-1$, (\ref{eq1-osw})-(\ref{eq2-osw}) is the Cordoba-Cordoba-Fontelos model introduced in \cite{CCF1} for the 2D quasi-geostrophic equation.  Cordoba, Cordoba and Fontelos \cite{CCF1, CCF2} showed finite-time singularity formation for this model with a general class of initial data. 

For axisymmetric 3D incompressible Navier-Stokes equation with swirl, Hou, Li, Shi, Wang and Yu \cite{Hou} proposed a 1D nonlocal model for a simplified 3D nonlocal system \cite{HouL}. For this 1D model, the authors proved finite-time singularity formation rigorously and showed numerical evidences.

\medskip

\subsection{1D models for MHD} 
% The toy models for the Hall MHD will be studied in a separate paper.
Inspired by the works discussed above, we will propose a family of nonlocal nonlinear models for the MHD system (\ref{mhd}) as an attempt to understand the intricate structures involved in this system. In the context of MHD, besides the convection and stretching effects, the coupling and interaction between the fluid velocity and magnetic field also play crucial roles, which naturally introduce additional challenges. 

Denote the Els\"asser variables by
\[p=u+B, \ \ m=u-B.\]
Equivalent to (\ref{mhd}), $(p,m)$ satisfies the system
\begin{equation}\label{mhd-el}
\begin{split}
p_t+(m\cdot\nabla) p+\nabla \Pi=&\ 0, \\
m_t+(p\cdot\nabla) m +\nabla \Pi=&\ 0, \\
\nabla\cdot p=0, \ \ \nabla\cdot m=&\ 0.
\end{split}
\end{equation}
The structure of system (\ref{mhd-el}) indicates that $p$ and $m$ are transported by each other. We also note that (\ref{mhd-el}) appears in a rather symmetric form. Denote the vorticity of $p$ and $m$ by
\[
\Omega=\nabla\times p, \ \ \omega=\nabla\times m.
\]
It follows from the Biot-Savart law that 
\begin{equation}\notag%\label{BS}
p=\nabla\times (-\Delta)^{-1} \Omega, \ \ m=\nabla\times (-\Delta)^{-1} \omega.
\end{equation}
Taking the curl $\nabla\times$ on the equations of (\ref{mhd-el}) gives
\begin{equation}\label{mhd-omega}
\begin{split}
\Omega_t+(m\cdot\nabla) \Omega -(\Omega\cdot\nabla) m+\nabla\times (m\nabla p)=&\ 0, \\
\omega_t +(p\cdot\nabla) \omega -(\omega\cdot\nabla) p+\nabla\times (p\nabla m)=&\ 0, \\
p=\nabla\times (-\Delta)^{-1} \Omega, \\
m=\nabla\times (-\Delta)^{-1} \omega,
\end{split}
\end{equation}
where $m\nabla p=(m_i\partial_j p_i)_j$ and $p\nabla m=(p_i\partial_j m_i)_j$.
We propose the following 1D model to mimic system (\ref{mhd-omega}), 
%on $[-\pi, \pi]$,
\begin{equation}\label{mhd-1d0}
\begin{split}
\Omega_t+m \Omega_x-\Omega m_x-\omega p_x+m\Omega_x=&\ 0, \\
\omega_t+p\omega_x-\omega p_x - \Omega m_x+p\omega_x=&\ 0, \\
p_x=H \Omega, \ \ 
m_x=&\ H \omega.
\end{split}
\end{equation}
In this paper, we will work with a simplified version of (\ref{mhd-1d0}) by dropping the stretching effects $\Omega m_x$ and $\omega p_x$ and focusing on the transport effects, namely
\begin{equation}\label{mhd-1d}
\begin{split}
\Omega_t+am \Omega_x-\omega p_x=&\ 0, \\
\omega_t+ap\omega_x - \Omega m_x=&\ 0, \\
p_x=H \Omega, \ \ 
m_x=&\ H \omega,
\end{split}
\end{equation}
with a parameter $a\in \mathbb R$.
%As in some historical study for the Euler equation, see \cite{}, the constitution relations $p_x=H \Omega$ and $m_x= H \omega$ are introduced to mimic the Biot-Savart laws (\ref{BS}) which are only valid in 2D and 3D. In (\ref{mhd-1d}), the magnitude of the parameter $a\in \mathbb R$ represents the strength of the transport effect and its sign tells the transport direction. 
We will investigate (\ref{mhd-1d}) on the periodic interval $S^1=[-\pi, \pi]$. Correspondingly, the Hilbert transform for periodic functions on $S^1$ can be defined as
\begin{equation}\label{H2}
H f(x)=\frac{1}{2\pi} P.V. \int_{-\pi}^{\pi}f(y)\cot\left( \frac{x-y}{2}\right)\, dy.
\end{equation}
Indeed, the Cauchy kernel $\frac{1}{x}$ in definition (\ref{H1}) can be made periodic using the following identity
\begin{equation}\notag
\frac12\cot\left( \frac{x-y}{2}\right)=\frac{1}{x}+\sum_{n=1}^{\infty}\left(\frac{1}{x+2n\pi}+\frac{1}{x-2n\pi}\right).
\end{equation}
To uniquely determine $p$ from $\Omega$ and $m$ from $\omega$, we make the choice of Gauge by taking zero-mean value
\begin{equation}\label{gauge1}
\int_{-\pi}^{\pi} p(t,x)\,dx =\int_{-\pi}^{\pi} m(t,x)\,dx=0.
\end{equation}
%or fixing one point
%\begin{equation}\label{gauge2}
%p(t,0) = m(t, 0)=0.
%\end{equation}
%{\color{red} }

%Notice that 
%\[u=\frac12(p+m), \ \ \ B=\frac12(p-m),\]
%and hence
%\begin{equation}\label{ub-pm}
%\begin{split}
%u_x=\frac12(p_x+m_x)=& \frac12 H(\Omega+\omega), \\
%B_x=\frac12(p_x-m_x)=& \frac12 H(\Omega-\omega).
%\end{split}
%\end{equation}

We note that the mean value of $\Omega$ and $\omega$ is invariant for system (\ref{mhd-1d}) with $a=1$. Indeed, we have for a smooth solution $(\Omega, \omega)$ that
\begin{equation}\notag
\begin{split}
\frac{d}{dt} \int_{-\pi}^{\pi}\Omega(t,x)\, dx=&\int_{-\pi}^{\pi}\left(-am\Omega_x+\omega p_x\right)\,dx\\
=&\int_{-\pi}^{\pi}\left(am_x\Omega+\omega p_x\right)\,dx\\
=&\int_{-\pi}^{\pi}\left(a\Omega H\omega+\omega H\Omega\right)\,dx\\
%=&\int_{-\pi}^{\pi}\left(-a\omega H\Omega+\omega H\Omega\right)\,dx\\
=&(1-a)\int_{-\pi}^{\pi}\omega H\Omega\,dx\\
\end{split}
\end{equation}
where we have used integration by parts and the skew symmetry property of the Hilbert transform. Similarly, we have
\begin{equation}\notag
\begin{split}
\frac{d}{dt} \int_{-\pi}^{\pi}\omega(t,x)\, dx=\int_{-\pi}^{\pi}\left(-ap\omega_x+\Omega m_x\right)\,dx
=(1-a)\int_{-\pi}^{\pi}\Omega H\omega\,dx.
\end{split}
\end{equation}
Obviously when $a=1$, it follows
\[\frac{d}{dt} \int_{-\pi}^{\pi}\Omega(t,x)\, dx=\frac{d}{dt} \int_{-\pi}^{\pi}\omega(t,x)\, dx=0,\]
and this is not true in general for $a\neq1$.
Hence, it is not appropriate to consider solutions of (\ref{mhd-1d}) in spaces of functions with zero mean for general value of $a$.

Consider the rescaled variables
\[\tilde \Omega=a\Omega, \ \ \ \tilde \omega=a\omega \]
with corresponding $\tilde p$ and $\tilde m$ such that
\[\tilde p_x=H\tilde \Omega, \ \ \tilde m_x=H\tilde \omega.\]
We can verify that $\tilde p=ap$ and $\tilde m=am$.
In view of (\ref{mhd-1d}), $(\tilde \Omega, \tilde \omega)$ satisfies the system
\begin{equation}\label{mhd-scale}
\begin{split}
\tilde \Omega_t+\tilde m\tilde\Omega_x-a^{-1}\tilde \omega\tilde p_x=&\ 0,\\
\tilde \omega_t+\tilde p\tilde\omega_x-a^{-1}\tilde \Omega\tilde m_x=&\ 0.
\end{split}
\end{equation}
Formally, taking $a\to\infty$, (\ref{mhd-scale}) turns to the system with only convection effect (with the tilde sign suppressed),
\begin{equation}\label{mhd-tran}
\begin{split}
 \Omega_t+ m\Omega_x=&\ 0,\\
 \omega_t+ p\omega_x=&\ 0, \\
 p_x=H\Omega, \ \ m_x=&\ H\omega.
\end{split}
\end{equation}

We will investigate both systems (\ref{mhd-1d}) and (\ref{mhd-tran}) in the paper. We point out that formulating the problem in Els\"asser variables does not give us essential advantage; rather it has the benefit of dealing with less nonlinear terms.  

%In this paper, we will study (\ref{mhd-1d}) and (\ref{mhd-tran}) on the periodic domain $[-\pi, \pi]$. 
%In view of (\ref{H2}), $p_x=H\Omega$ and $m_x=H\omega$, we have
%\begin{equation}\label{eq-pm}
%\begin{split}
%p(t,x)=&\ \frac{1}{\pi} P. V. \int_{-\pi}^{\pi}\Omega(t,y)\log \left|\sin\left(\frac{x-y}{2}\right)\right|\, dy, \\
%m(t,x)=&\ \frac{1}{\pi} P. V. \int_{-\pi}^{\pi}\omega(t,y)\log \left|\sin\left(\frac{x-y}{2}\right)\right|\, dy.
%\end{split}
%\end{equation} 

\medskip

\subsection{Main results}
For general $a\in\mathbb R$, we show the existence of local in time solutions to (\ref{mhd-1d}) in the space $\mathcal H^1(S^1)$.
\begin{Theorem}\label{thm-existence}
Let $a\in \mathbb R$ and $\Omega_0, \omega_0\in \mathcal H^1(S^1)$. There exists a time $T$ which depends on $\|\Omega_{0,x}\|_{L^2}$ and $\|\omega_{0,x}\|_{L^2}$ such that there exists a unique solution $(\Omega(t,x),$ $ \omega(t,x))$ to (\ref{mhd-1d}) with initial data $\Omega(0,x)=\Omega_0$ and $\omega(0,x)=\omega_0$ on $[0,T)$, which satisfies
\[\Omega, \omega\in C^0\left([0,T); \mathcal H^1(S^1)\right)\cap C^1\left([0,T); L^2(S^1)\right).\]
\end{Theorem}

The following theorem provides a Beale-Kato-Majda type of regularity criterion.

\begin{Theorem}\label{thm-criterion}
Let $(\Omega(t,x), \omega(t,x))$ be the solution of (\ref{mhd-1d}) on $[0,T)$ obtained in Theorem \ref{thm-existence}. If 
\begin{equation}\label{criterion}
\int_0^T\left(\|H\Omega (t)\|_{L^\infty}+\|H\omega(t)\|_{L^\infty}\right) \, dt<\infty,
\end{equation}
the solution can be extended beyond $T$ in the space $\mathcal H^1(S^1)\times \mathcal H^1(S^1)$.
\end{Theorem}

Furthermore, if the initial data is in a space with higher regularity, the solution obtained in Theorem \ref{thm-existence} also has higher regularity. Specifically, we will show:
\begin{Theorem}\label{thm-higher}
Assume $\Omega_0, \omega_0\in \mathcal H^n(S^1)$ with $n\geq 2$. Let $(\Omega, \omega)$ be a solution of (\ref{mhd-1d}) with initial data $(\Omega_0, \omega_0)$ on $[0,T)$, satisfying $\Omega, \omega\in C([0,T); \mathcal H^1)$. Then, we have
\begin{equation}\notag
\sup_{0\leq t< T}\left(\|\Omega(t)\|_{\mathcal H^n}+\|\omega(t)\|_{\mathcal H^n}\right)<\infty.
\end{equation}
\end{Theorem}

With the absence of stretching effect, the solution of (\ref{mhd-tran}) can be shown to exist in the space $H^1(S^1)$ for all the time. Namely, we have 
\begin{Theorem}\label{thm-tran}
Assume $\Omega_0,\omega_0\in \mathcal H^1(S^1)$. Then there exists a unique solution $(\Omega(t),\omega(t))$ of (\ref{mhd-tran}) with initial data $(\Omega_0,\omega_0)$ on $[0,\infty)$.
\end{Theorem}

%\begin{Remark}\label{rk-numerics}
%Numerical evidences of singularity and global existence: (i) small $a$, singularity; (ii) large $a$, global existence.
%\end{Remark}
Some numerical simulations will be provided in Section \ref{sec-numerics}. The numerical results suggest that starting from smooth periodic initial data, solutions of the model (\ref{mhd-1d}) with $a=1$ or $a=-1$ are unlikely to develop singularities at finite time. This observation agrees with the numerical results done by De Gregorio \cite{DeG1, DeG2} and Okamoto, Sakajo and Wunsch \cite{OSW} for the De Gregorio model (\ref{eq1-deg})-(\ref{eq2-deg}).

\bigskip

\section{Notations and preliminaries}

\subsection{Functional setting}
Denote 
\begin{equation}\notag
\begin{split}
L^2(S^1)=& \left\{f| f\in L^2(-\pi, \pi), f \ \mbox{is periodic on} [-\pi,\pi]\right\},\\
\mathcal H^k(S^1)=& \left\{f| f^{(s)}\in L^2(-\pi, \pi),  f^{(s)} \ \mbox{is periodic on} [-\pi,\pi], \ \ \mbox{for all} \ \ 0\leq s\leq k\right\}.
\end{split}
\end{equation}
In particular, we consider the triplet of spaces
\[\mathcal V=\left\{ f\middle| f\in \mathcal H^2(S^1), \ \int_{-\pi}^{\pi} f(x)\, dx=0 \right\}\]
\[\mathcal W=\mathcal H^1(S^1), \ \ \mathcal X= L^2(S^1), \ \ \]
with the obvious embedding $\mathcal V\subset \mathcal W\subset \mathcal X$.

We denote $(,)$ by 
\[(f, g)=\int^{\pi}_{-\pi} fg\, dx.\]
The space $\mathcal H^k(S^1)$ is a Hilbert space endowed with the natural inner product 
\begin{equation}\notag
(f, g)_{\mathcal H^k}=\sum_{s=0}^k\left(f^{(s)}, g^{(s)}\right) \ \ \mbox{ for functions} \ \ f, g\in \mathcal H^k(S^1),
\end{equation}
and norm $(f, f)_{\mathcal H^k}^{\frac12}$.
%To specify the space, the notation $(,)_{\mathcal Z}$ is also used for the inner product of the Hilbert space $\mathcal Z$.

A bilinear form $\left<,\right>: \mathcal V\times \mathcal X\to \mathbb R$ is defined as
\[\left<f, g\right>=-\int^{\pi}_{-\pi} f_{xx}g\, dx.\]
Applying the integration by parts, we have for all $f\in \mathcal V$ and $g\in \mathcal W$
\[\left<f, g\right>=(f_x, g_x).\]

For a space $\mathcal Z$, we denote $\mathcal Z^2=\mathcal Z\times \mathcal Z$ by convention. In the context of a coupled system, for instance (\ref{mhd-1d}), it is convenient to introduce the triplet $\{\mathcal V^2, \mathcal W^2, \mathcal X^2\}$. Naturally, the Hilbert space $\mathcal W^2$ is endowed with the inner product 
\begin{equation}\notag
(f, g)_{\mathcal W^2}=(f_{1}, g_{1})_{\mathcal W} +(f_{2}, g_{2})_{\mathcal W} \ \ \ \forall \ \ f=(f_1, f_2)\in \mathcal W^2, \ \ g=(g_1, g_2)\in \mathcal W^2.
\end{equation}
%and the norm $(f_x, f_x)^{\frac12}$. 
In an analogous way, inner product can be defined for $\mathcal V^2$ and $\mathcal X^2$. A bilinear form $\left<,\right>:\mathcal V^2\times \mathcal X^2\to \mathbb R$ is defined as
\begin{equation}\label{bilin}
\left<f, g\right>=-\int^{\pi}_{-\pi} f_{1, xx}g_1\, dx-\int^{\pi}_{-\pi} f_{2, xx}g_2\, dx.
\end{equation}
For all $f=(f_1, f_2)\in\mathcal V^2$ and $g=(g_1, g_2)\in \mathcal W^2$,
we also have 
\[\left<f, g\right>=(f_{1,x}, g_{1,x})+(f_{2,x}, g_{2,x}).\]

\begin{Definition}\label{def-adm}
A family $\{\mathcal Z, \mathcal H, \mathcal Y\}$ of three real separable Banach spaces is called an admissible triplet if the following conditions hold:\\
(i) The inclusions $\mathcal Z\subset \mathcal H\subset \mathcal Y$ are continuous and dense.\\
(ii) $H$ is a Hilbert space endowed with inner product $(,)_{\mathcal H}$ and norm $\|\|_{\mathcal H}=(,)_{\mathcal H}^{\frac12}$.\\
(iii) There is a continuous non-degenerate bilinear form on $\mathcal Z\times \mathcal Y$, denoted by $\left<,\right>$, such that 
\begin{equation}\label{cond3}
\left<v, u\right>=(v, u)_{\mathcal H}, \ \ \mbox{for} \ \ v\in \mathcal V \ \ \mbox{and} \ \ u\in \mathcal H.
\end{equation}
\end{Definition} 

Denote $C_w$ by the space of functions with weak continuity and $C^1_w$ the space of functions with weak differentiability. 

An abstract theorem of existence of Kato-Lai \cite{KL} is stated as follows.
\begin{Theorem}\label{thm-KL}
Let $\{\mathcal Z, \mathcal H, \mathcal Y\}$ be an admissible triplet. Let $A: \mathcal H\to \mathcal Y$ be a weakly continuous map such that
\begin{equation}\label{cond-A}
\left<v, A(v)\right>\geq -\beta (\|v\|_{\mathcal H}^2) , \ \ \forall \ \ v\in \mathcal Z
\end{equation}
where $\beta(r)\geq 0$ is a monotone increasing function of $r\geq 0$. Then for any $u_0\in \mathcal H$, there exists a time $T>0$ such that the Cauchy problem
\begin{equation}\notag
u_t+A(u)=0, \ \ \ u(0,x)=u_0
\end{equation}
has a solution $u(t,x)$ on $[0,T]$ satisfying
\[u\in C_w([0,T]; \mathcal H)\cap C^1_w([0,T]; \mathcal Y).\]
Moreover, $\sup_{0<t<T}\|u(t)\|_H$ depends only on $T$, $\beta$ and $\|u_0\|_{\mathcal H}$.
\end{Theorem}

In order to prove the existence part of Theorem \ref{thm-existence}, the Kato-Lai theorem will be applied to system (\ref{mhd-1d}) with the admissible triplet $\{\mathcal V^2, \mathcal W^2, \mathcal X^2\}$.

\subsection{Properties of Hilbert transform}
The Hilbert transform has the following simple properties
\begin{equation}\notag
\begin{split}
H(cf)=&\ cHf, \ \ \mbox{for a constant} \ \ c,\\
H\sin(kx)=&-\cos(kx),\ \ \ H\cos(kx)=\sin(kx).
\end{split}
\end{equation}
And more generally, we have
\begin{equation}\notag
H\sin(kx+\theta)=-\cos(kx+\theta), \ \ H\cos(kx+\theta)=\sin(kx+\theta).
\end{equation}

For any periodic function $f$, the mean value of its Hilbert transform is zero, that is
\begin{equation}\label{Hm}
\int_{-\pi}^{\pi}Hf\, dx=0.
\end{equation}

\begin{Lemma}\cite{Zy}
The Hilbert transform $H$ is a bounded linear operator from space $L^p$ to $L^p$ with $1<p<\infty$ and
\begin{equation}\label{HLp}
\|Hf\|_{L^p}\leq C_p \|f\|_{L^p}
\end{equation}
for a constant $C_p>0$ depending on $p$.
\end{Lemma}

\bigskip

\section{Local existence}
\label{sec-existence}

This section is devoted to a proof of Theorem \ref{thm-existence}.  The proof includes three steps: (i) establishing the local existence of a solution by employing Theorem \ref{thm-KL}; (ii) showing the uniqueness of solution by a rather standard argument; (iii) justifying the strong continuity which is a consequence of the uniqueness and the time-reversible property of system (\ref{mhd-1d}).

{\textbf{Proof of Theorem \ref{thm-existence}:}}
Denote $u=(\Omega, \omega)$, $q=(p,m)$, and naturally $q_x=Hu=(H\Omega, H\omega)$. Denote $A(u)=(A_1(u), A_2(u))$ with 
\[A_1(u)= am\Omega_x-\omega p_x, \ \ \ A_2(u)= ap\omega_x-\Omega m_x.\]
Thus, system (\ref{mhd-1d}) can be written as
\begin{equation}\notag%\label{mhd-ab}
u_t+A(u)=0.
\end{equation}
It is obvious that the family $\{\mathcal V^2, \mathcal W^2, \mathcal X^2\}$ is an admissible triplet associated with the bilinear form $\left<,\right>$ defined in (\ref{bilin}). 
To apply Theorem \ref{thm-KL}, we will need to show that the operator $A$ maps $\mathcal W^2$ into $\mathcal X^2$ continuously and it satisfies (\ref{cond-A}). Indeed, for any $u=(\Omega, \omega)\in \mathcal W^2$ with $q=(p,m)\in \mathcal V^2$, we have
\begin{equation}\notag
\begin{split}
\|A(u)\|_{\mathcal X^2}=& \left(\|am\Omega_x-\omega p_x\|_{L^2}^2+\|ap\omega_x-\Omega m_x\|_{L^2}^2\right)^{\frac12}\\
\leq &\ \|am\Omega_x-\omega p_x\|_{L^2}+\|ap\omega_x-\Omega m_x\|_{L^2}\\
\leq &\ |a|\|m\|_{L^\infty} \|\Omega_x\|_{L^2}+\|\omega\|_{L^\infty} \|p_x\|_{L^2}\\
&+|a|\|p\|_{L^\infty} \|\omega_x\|_{L^2}+\|\Omega\|_{L^\infty} \|m_x\|_{L^2}\\
\leq &\ c_0\left(|a|\|m_x\|_{L^2} \|\Omega_x\|_{L^2}+\|\omega\|_{\mathcal H^1} \|p_x\|_{L^2}\right.\\
&\left.+|a|\|p_x\|_{L^2} \|\omega_x\|_{L^2}+\|\Omega\|_{\mathcal H^1} \|m_x\|_{L^2}\right)\\
\leq &\ c_0\left(|a|+1\right)\left(\|H\omega\|_{L^2} \|\Omega\|_{\mathcal H^1}+\|\omega\|_{\mathcal H^1} \|H\Omega\|_{L^2}\right)\\
\leq &\ c_0\left(|a|+1\right)\left(\|\omega\|_{L^2} \|\Omega\|_{\mathcal H^1}+\|\omega\|_{\mathcal H^1} \|\Omega\|_{L^2}\right)\\
\end{split}
\end{equation}
where we have used the H\"older inequality, Sobolev inequality, the fact that $p$ and $m$ have zero mean, and the property (\ref{HLp}). 
% the Hilbert transform $H$ is a bounded operator from space $L^p$ to $L^p$. 
It follows that $A$ maps $\mathcal W^2$ into $\mathcal X^2$. On the other hand, for any $u_1=(\Omega_1, \omega_1)\in \mathcal W^2$ with $q_1=(p_1,m_1)\in \mathcal V^2$ and $u_2=(\Omega_2, \omega_2)\in \mathcal W^2$ with $q_2=(p_2,m_2)\in \mathcal V^2$, we deduce 
\begin{equation}\label{est-u1u2-1}
\begin{split}
\|A(u_1)-A(u_2)\|_{\mathcal X^2}=&\left(\|(am_1\Omega_{1,x}-\omega_1 p_{1,x})-(am_2\Omega_{2,x}-\omega_2 p_{2,x})\|_{L^2}^2\right.\\
&\left. + \|(ap_1\omega_{1,x}-\Omega_1 m_{1,x})-(ap_2\omega_{2,x}-\Omega_2 m_{2,x})\|_{L^2}^2\right)^{\frac12}\\
\leq &\ \|(am_1\Omega_{1,x}-\omega_1 p_{1,x})-(am_2\Omega_{2,x}-\omega_2 p_{2,x})\|_{L^2}\\
&+  \|(ap_1\omega_{1,x}-\Omega_1 m_{1,x})-(ap_2\omega_{2,x}-\Omega_2 m_{2,x})\|_{L^2}.
\end{split}
\end{equation}
Applying the H\"older inequality, Sobolev inequality, and (\ref{HLp}) leads to
\begin{equation}\label{est-u1u2-2}
\begin{split}
&\|(am_1\Omega_{1,x}-\omega_1 p_{1,x})-(am_2\Omega_{2,x}-\omega_2 p_{2,x})\|_{L^2}\\
\leq& \ |a|\|\Omega_{1,x}\|_{L^2}\|m_1-m_2\|_{L^\infty}+|a|\|\Omega_{1,x}-\Omega_{2,x}\|_{L^2}\|m_2\|_{L^\infty}\\
&+\|\omega_{2}\|_{L^\infty}\|p_{2,x}-p_{1,x}\|_{L^2}+|a|\|\omega_{2}-\omega_{1}\|_{L^\infty}\|p_{1,x}\|_{L^2}\\
\leq& \ c_0|a|\|\Omega_{1,x}\|_{L^2}\|m_{1,x}-m_{2,x}\|_{L^2}+c_0|a|\|\Omega_{1,x}-\Omega_{2,x}\|_{L^2}\|m_{2,x}\|_{L^2}\\
&+c_0\|\omega_{2}\|_{\mathcal H^1}\|p_{2,x}-p_{1,x}\|_{L^2}+c_0|a|\|\omega_{2}-\omega_{1}\|_{\mathcal H^1}\|p_{1,x}\|_{L^2}\\
\leq &\ c_0(|a|+1) \left(\|\Omega_1\|_{\mathcal H^1}+\|\omega_2\|_{\mathcal H^1}\right)
 \left(\|\Omega_1-\Omega_2\|_{\mathcal H^1}+\|\omega_1-\omega_2\|_{\mathcal H^1}\right),
%\leq &\ c_0(|a|+1) \left(\|\Omega_1\|_{L^2}+\|\omega_2\|_{L^2}+\|\Omega_{1,x}\|_{L^2}+\|\omega_{2,x}\|_{L^2}\right)\\
%&\cdot \left(\|\Omega_1-\Omega_2\|_{L^2}+\|\omega_1-\omega_2\|_{L^2}+\|\Omega_{1,x}-\Omega_{2,x}\|_{L^2}+\|\omega_{1,x}-\omega_{2,x}\|_{L^2}\right),
\end{split}
\end{equation}
and similarly 
\begin{equation}\label{est-u1u2-3}
\begin{split}
&\|(ap_1\omega_{1,x}-\Omega_1 m_{1,x})-(ap_2\omega_{2,x}-\Omega_2 m_{2,x})\|_{L^2}\\
\leq &\ c_0(|a|+1)  \left(\|\Omega_1\|_{\mathcal H^1}+\|\omega_2\|_{\mathcal H^1}\right)
 \left(\|\Omega_1-\Omega_2\|_{\mathcal H^1}+\|\omega_1-\omega_2\|_{\mathcal H^1}\right).
%\leq &\ c_0(|a|+1) \left(\|\Omega_1\|_{L^2}+\|\omega_1\|_{L^2}+\|\Omega_{2,x}\|_{L^2}+\|\omega_{2,x}\|_{L^2}\right)\\
%&\cdot \left(\|\Omega_1-\Omega_2\|_{L^2}+\|\omega_1-\omega_2\|_{L^2}+\|\Omega_{1,x}-\Omega_{2,x}\|_{L^2}+\|\omega_{1,x}-\omega_{2,x}\|_{L^2}\right).
\end{split}
\end{equation}
The estimates (\ref{est-u1u2-1})-(\ref{est-u1u2-3}) together indicate that $A: \mathcal W^2\to \mathcal X^2$ is strongly continuous. 

By the definition of the bilinear form in (\ref{bilin}), we have for any $u=(\Omega, \omega)\in \mathcal V^2$
\begin{equation}\label{est-uau1}
\begin{split}
\left<u, A(u)\right>=&- (u_{xx}, A(u))=(u_x, (A(u))_x)\\
=&\ (\Omega_x, (am\Omega_x-\omega p_x)_x)+(\omega_x, (ap\omega_x-\Omega m_x)_x)\\
=& \int_{-\pi}^{\pi} \Omega_x (am_x\Omega_x+am\Omega_{xx}-\omega_x p_x-\omega p_{xx})\, dx\\
&+\int_{-\pi}^{\pi} \omega_x (ap_x\omega_x+ap\omega_{xx}-\Omega_x m_x-\Omega m_{xx})\, dx.
\end{split}
\end{equation} 
Note that $A(u)\in \mathcal X^2$ and (\ref{est-uau1}) can be made rigorous through a standard approximating procedure. 
Applying integration by parts to the right hand side of (\ref{est-uau1}), it has
\begin{equation}\notag
a\int_{-\pi}^{\pi} m \Omega_x \Omega_{xx} \,dx=-a\int_{-\pi}^{\pi} m_x \Omega_x \Omega_{x} \,dx-a\int_{-\pi}^{\pi} m \Omega_{xx} \Omega_{x} \,dx.
\end{equation}
Hence we conclude 
\begin{equation}\label{est-uau2}
a\int_{-\pi}^{\pi} m \Omega_x \Omega_{xx} \,dx=-\frac{a}2\int_{-\pi}^{\pi} m_x \Omega_x^2\,dx,
\end{equation}
and similarly
\begin{equation}\label{est-uau3}
a\int_{-\pi}^{\pi} p \omega_x \omega_{xx} \,dx=-\frac{a}2\int_{-\pi}^{\pi} p_x \omega_x^2\,dx.
\end{equation}
Since $p_x=H\Omega$ and $m_x=H\omega$, combining (\ref{est-uau1})-(\ref{est-uau3}) gives
\begin{equation}\label{est-uau4}
\begin{split}
\left<u, A(u)\right>=&\ (u_x, (A(u))_x)\\
=&\ \frac{a}2\int_{-\pi}^{\pi} (H\omega) \Omega_x^2\,dx+\frac{a}2\int_{-\pi}^{\pi} (H\Omega) \omega_x^2\,dx\\
&-\int_{-\pi}^{\pi}  \Omega_x\omega_x H\Omega\,dx-\int_{-\pi}^{\pi}  \omega\Omega_x H\Omega_x\,dx\\
&-\int_{-\pi}^{\pi}  \Omega_x\omega_x H\omega\,dx-\int_{-\pi}^{\pi}  \Omega\omega_x H\omega_x\,dx.
\end{split}
\end{equation} 
Applying H\"older's inequality, Sobolev's inequality, (\ref{Hm}) and (\ref{HLp}), we have
\begin{equation}\label{est-uau5}
\begin{split}
\left|\int_{-\pi}^{\pi} H\omega \Omega_x^2\,dx\right|\leq&\ \|H\omega\|_{L^\infty} \|\Omega_x\|_{L^2}^2\\
\leq&\ c_0 \|H\omega_x\|_{L^2} \|\Omega_x\|_{L^2}^2\\
\leq&\ c_0 \|\omega_x\|_{L^2} \|\Omega_x\|_{L^2}^2,
\end{split}
\end{equation}
and similarly 
\begin{equation}\label{est-uau6}
\begin{split}
\left|\int_{-\pi}^{\pi}\Omega_x\omega_x H\Omega\,dx\right|+\left|\int_{-\pi}^{\pi}\omega\Omega_x H\Omega_x\,dx\right|\leq&\ c_0 \|\omega\|_{\mathcal H^1} \|\Omega_x\|_{L^2}^2,\\
\left|\int_{-\pi}^{\pi}(H\Omega) \omega_x^2\,dx\right| \leq&\ c_0 \|\omega_x\|_{L^2}^2 \|\Omega_x\|_{L^2},\\
\left|\int_{-\pi}^{\pi}\omega_x\Omega_x H\omega\,dx\right|
+\left|\int_{-\pi}^{\pi}\omega_x\Omega H\omega_x\,dx\right|
\leq&\ c_0 \|\omega_x\|_{L^2}^2 \|\Omega\|_{\mathcal H^1}.
\end{split}
\end{equation}
Therefore, putting together (\ref{est-uau4})-(\ref{est-uau6}), we deduce
\begin{equation}\label{est-uau7}
\begin{split}
\left|\left<u, A(u)\right>\right|\leq &\ c_0(|a|+1) \left(\|\omega\|_{\mathcal H^1} \|\Omega_x\|_{L^2}^2+\|\omega_x\|_{L^2}^2 \|\Omega\|_{\mathcal H^1}\right)\\
\leq&\ c_0(|a|+1) \left(\|\omega\|_{\mathcal H^1} + \|\Omega\|_{\mathcal H^1}\right)^3. 
\end{split}
\end{equation}
Hence, the operator $A$ satisfies (\ref{cond-A}) with $\beta(r)=c_0(|a|+1) r^{\frac32}$. As a consequence, applying Theorem \ref{thm-KL}, we conclude that there exists a time $T>0$ such that system (\ref{mhd-1d}) has a solution $\left(\Omega(t,x), \omega(t,x)\right)$ on $[0,T]$ satisfying 
\[\Omega, \omega\in C_w([0,T]; \mathcal W)\cap C^1_w([0,T]; \mathcal X).\] 

Next we show the uniqueness of solution to (\ref{mhd-1d}). Let $u_1=(\Omega_1, \omega_1)$ be a solution to (\ref{mhd-1d}) with initial data $u_0=(\Omega_0, \omega_0)$. Let $q_1=(p_1, m_1)$ such that $p_{1,x}=H\Omega_1$
and $m_{1,x}=H\omega_1$.  Let $u_2=(\Omega_2, \omega_2)$ be another solution to (\ref{mhd-1d}) with the same initial data $(\Omega_0, \omega_0)$ and associated with $q_2=(p_2, m_2)$. Since both $(\Omega_1, \omega_1)$
and $(\Omega_2, \omega_2)$ satisfy (\ref{mhd-1d}), we are able to show that (details omitted)
\begin{equation}\label{est-diff1}
\begin{split}
&\frac12\frac{d}{dt}\left(\|\Omega_1(t)-\Omega_2(t)\|_{L^2}^2+\|\omega_1(t)-\omega_2(t)\|_{L^2}^2\right)\\
\leq&c_0 (|a|+1) \max_{0\leq t\leq T}\left(\|\Omega_{1}\|_{\mathcal H^1}+\|\Omega_{2}\|_{\mathcal H^1}+\|\omega_{1}\|_{\mathcal H^1}+\|\omega_{2}\|_{\mathcal H^1}\right)\\
& \cdot \left(\|\Omega_1(t)-\Omega_2(t)\|_{L^2}^2+\|\omega_1(t)-\omega_2(t)\|_{L^2}^2\right).
\end{split}
\end{equation}
Thus, uniqueness follows from (\ref{est-diff1}) and Gr\"onwall's inequality.

Strong continuity in time follows from the uniqueness and the fact that system (\ref{mhd-1d}) is time-reversible. Indeed,   %it is clear that $\Omega$ and $\omega$ are strongly continuously at $t=0$. 
it follows from (\ref{est-uau7}) that
 \[\|\Omega_x(t)\|_{L^2}+\|\omega_x(t)\|_{L^2}\to \|\Omega_{0,x}\|_{L^2}+\|\omega_{0,x}\|_{L^2} \ \ \mbox{as} \ \ t\to 0.\]
Hence, we know
\[\Omega(t)\to \Omega_0, \ \ \omega(t)\to \omega_0 \ \ \mbox{strongly in} \ \ \mathcal H^1 \ \ \mbox{as} \ \ t\to 0.\]
As a consequence of uniqueness, $\Omega$ and $\omega$ are strongly right-continuous. In addition, the property of time-reversibility implies that $\Omega$ and $\omega$ are strongly left-continuous as well.

\cbdu

\bigskip

\section{Regularity criterion}
\label{sec-criterion}

In this section, we prove Theorem \ref{thm-criterion} and the higher regularity result in Theorem \ref{thm-higher}.

{\textbf{Proof of Theorem \ref{thm-criterion}:} In view of the local existence theorem, we just need to show that the $\mathcal H^1$ norm of $\Omega(t)$ and $\omega(t)$ remains bounded as $t\to T$ under condition (\ref{criterion}).  

Assume $u=(\Omega, \omega)$ is a solution of (\ref{mhd-1d}) on $[0,T)$.  We note that
\begin{equation}\label{est-h11}
\begin{split}
&\frac12\frac{d}{dt}\left(\|\Omega_x\|_{L^2}^2+\|\omega_x\|_{L^2}^2\right)\\
=&\left(\Omega_x, \Omega_{tx}\right)+\left(\omega_x, \omega_{tx}\right)\\
=&\left(\Omega_x, -(am\Omega_x-\omega p_x)_x\right)+\left(\omega_x, -(ap\omega_x-\Omega m_x)_x\right)\\
=& -\frac{a}2\int_{-\pi}^{\pi} (H\omega) \Omega_x^2\,dx-\frac{a}2\int_{-\pi}^{\pi} (H\Omega) \omega_x^2\,dx\\
&+\int_{-\pi}^{\pi}  \Omega_x\omega_x H\Omega\,dx+\int_{-\pi}^{\pi}  \omega\Omega_x H\Omega_x\,dx\\
&+\int_{-\pi}^{\pi}  \Omega_x\omega_x H\omega\,dx+\int_{-\pi}^{\pi}  \Omega\omega_x H\omega_x\,dx
\end{split}
\end{equation}
where we used (\ref{est-uau4}) in the last step. Applying the identities 
\[(v, u)=(Hv, Hu), \ \ \ H(vHv)=\frac12\left((Hv)^2-v^2\right),\]
we infer
\begin{equation}\label{est-h12}
\begin{split}
\int_{-\pi}^{\pi}  \omega\Omega_x H\Omega_x\,dx=& \int_{-\pi}^{\pi}  (H\omega) H(\Omega_x H\Omega_x)\,dx\\
=&\ \frac12 \int_{-\pi}^{\pi}  (H\omega) \left( (H\Omega_x)^2-(\Omega_x)^2\right)\,dx,\\
\end{split}
\end{equation}
\begin{equation}\label{est-h13}
\begin{split}
\int_{-\pi}^{\pi}  \Omega\omega_x H\omega_x\,dx=& \int_{-\pi}^{\pi}  (H\Omega) H(\omega_x H\omega_x)\,dx\\
=&\ \frac12 \int_{-\pi}^{\pi}  (H\Omega) \left( (H\omega_x)^2-(\omega_x)^2\right)\,dx.\\
\end{split}
\end{equation}
Combining (\ref{est-h11})-(\ref{est-h13}), we have
\begin{equation}\label{est-h14}
\begin{split}
&\frac12\frac{d}{dt}\left(\|\Omega_x\|_{L^2}^2+\|\omega_x\|_{L^2}^2\right)\\
=& -\frac{a+1}2\int_{-\pi}^{\pi} (H\omega) \Omega_x^2\,dx-\frac{a+1}2\int_{-\pi}^{\pi} (H\Omega) \omega_x^2\,dx\\
&+\frac12\int_{-\pi}^{\pi} (H\omega) (H\Omega_x)^2\,dx+\frac12\int_{-\pi}^{\pi} (H\Omega) (H\omega_x)^2\,dx\\
&+\int_{-\pi}^{\pi}  \Omega_x\omega_x \left(H\Omega+H\omega\right)\,dx\\
\leq& \frac{|a+1|}{2} \|H\omega\|_{L^\infty}\|\Omega_x\|_{L^2}^2+\frac{|a+1|}{2} \|H\Omega\|_{L^\infty}\|\omega_x\|_{L^2}^2\\
&+ \frac12 \|H\omega\|_{L^\infty}\|\Omega_x\|_{L^2}^2+\frac{1}{2} \|H\Omega\|_{L^\infty}\|\omega_x\|_{L^2}^2\\
&+\|H\Omega+H\omega\|_{L^\infty}\|\Omega_x\|_{L^2}\|\omega_x\|_{L^2}\\
\leq &\ c_0(|a|+1) \left(\|H\Omega\|_{L^\infty}+\|H\omega\|_{L^\infty}\right)\left(\|\Omega_x\|_{L^2}^2+\|\omega_x\|_{L^2}^2\right)
\end{split}
\end{equation}
for a constant $c_0>0$. It follows from Gr\"onwall's inequality that 
\begin{equation}\notag
\begin{split}
&\left(\|\Omega_x(t)\|_{L^2}^2+\|\omega_x(t)\|_{L^2}^2\right)\\
\leq &\left(\|\Omega_x(0)\|_{L^2}^2+\|\omega_x(0)\|_{L^2}^2\right) \exp \left\{2c_0(|a|+1)\int_0^t\left(\|H\Omega (\tau)\|_{L^\infty}+\|H\omega(\tau)\|_{L^\infty}\right) \, d\tau \right\}.
\end{split}
\end{equation}
Thus, the statement of the theorem is justified.

\cbdu

{\textbf{Proof of Theorem \ref{thm-higher}:}} The statement can be established through standard energy method. We only deal with the case of $n=2$ and obtain the a priori estimate for $\|\Omega(t)\|_{\mathcal H^2}$ and $\|\omega(t)\|_{\mathcal H^2}$.
Formally, differentiating the equations of (\ref{mhd-1d}) twice in space yields
\begin{equation}\label{eq-xx}
\begin{split}
\Omega_{txx}=&-2am_x\Omega_{xx}+2\omega_xp_{xx}-am_{xx}\Omega_x\\
&+\omega_{xx}p_x-am\Omega_{xxx}+\omega p_{xxx},\\
\omega_{txx}=&-2a p_x\omega_{xx}+2\Omega_x m_{xx}-ap_{xx}\omega_x\\
&+\Omega_{xx} m_x-ap \omega_{xxx}+\Omega m_{xxx}.
\end{split}
\end{equation}
Taking the inner product of the first equation with $\Omega_{xx}$ and the second one with $\omega_{xx}$, we have
\begin{equation}\label{est-xx1}
\begin{split}
&\frac12\frac{d}{dt}\left(\|\Omega_{xx}\|_{L^2}^2+\|\omega_{xx}\|_{L^2}^2\right)\\
=&-2a(\Omega_{xx}, m_x\Omega_{xx})+2(\Omega_{xx},\omega_xp_{xx})-a(\Omega_{xx}, m_{xx}\Omega_x)\\
&+(\Omega_{xx}, \omega_{xx}p_x)-a(\Omega_{xx}, m\Omega_{xxx})+(\Omega_{xx}, \omega p_{xxx})\\
&-2a(\omega_{xx}, p_x\omega_{xx})+(\omega_{xx}, \Omega_x m_{xx})-a(\omega_{xx}, p_{xx}\omega_x)\\
&+(\omega_{xx},\Omega_{xx}m_x)-a(\omega_{xx}, p\omega_{xxx})+(\omega_{xx}, \Omega m_{xxx}).
\end{split}
\end{equation}
Notice that, by integration by parts,
\begin{equation}\notag
-a(\Omega_{xx}, m\Omega_{xxx})=a(\Omega_{xxx}, m\Omega_{xx})+a(\Omega_{xx}, m_x\Omega_{xx})
\end{equation}
which implies
\begin{equation}\notag
-a(\Omega_{xx}, m\Omega_{xxx})=\frac{a}{2}(\Omega_{xx}, m_x\Omega_{xx}).
\end{equation}
Similarly, we have
\begin{equation}\notag
-a(\omega_{xx}, p\omega_{xxx})=\frac{a}{2}(\omega_{xx}, p_x\omega_{xx}).
\end{equation}
Applying H\"older's inequality, the Hilbert transform boundedness on $L^p$, it follows
\begin{equation}\notag
\begin{split}
\left|(\Omega_{xx}, m_x\Omega_{xx})\right|=&\left|(\Omega_{xx}, (H\omega) \Omega_{xx})\right|\\
\leq &\ c_0 \|H\omega\|_{L^\infty}\|\Omega_{xx}\|_{L^2}^2\\
\leq &\ c_0 \|\omega_x\|_{L^2}\|\Omega_{xx}\|_{L^2}^2,\\
\end{split}
\end{equation}
and similarly
\begin{equation}\notag
\left|(\omega_{xx}, p_x\omega_{xx})\right|=\left|(\Omega_{xx}, (H\omega) \Omega_{xx})\right|
\leq c_0 \|\Omega_x\|_{L^2}\|\omega_{xx}\|_{L^2}^2.
\end{equation}
We estimate $(\Omega_{xx},\omega_xp_{xx})$ as
\begin{equation}\notag
\begin{split}
|(\Omega_{xx},\omega_xp_{xx})|=&\ |(\Omega_{xx},\omega_xH\Omega_{x})|\\
\leq &\ c_0 \|\Omega_{xx}\|_{L^2}\|\omega_x\|_{L^4}\|H\Omega_x\|_{L^4}\\
\leq &\ c_0 \|\Omega_{xx}\|_{L^2}\|\omega_x\|_{L^2}^{\frac34}\|\omega_{xx}\|_{L^2}^{\frac14}\|H\Omega_x\|_{L^2}^{\frac34}\|H\Omega_{xx}\|_{L^2}^{\frac14}\\
\leq&\ c_0 \|\Omega_{xx}\|_{L^2}^{\frac54}\|H\Omega_x\|_{L^2}^{\frac34}\|\omega_x\|_{L^2}^{\frac34}\|\omega_{xx}\|_{L^2}^{\frac14}\\
\leq& \ c_0  \|\Omega_{xx}\|_{L^2}^{\frac{15}{8}}\|H\Omega_x\|_{L^2}^{\frac98}+c_0 \|\omega_x\|_{L^2}^{\frac94}\|\omega_{xx}\|_{L^2}^{\frac34}\\
\leq&\ c_0 \|\Omega_x\|_{L^2}\|\Omega_{xx}\|_{L^2}^{2}+c_0 \|\omega_x\|_{L^2}\|\omega_{xx}\|_{L^2}^{2}, 
\end{split}
\end{equation}
where we used the inequalities of H\"older, Galiardo-Nirenberg and Young, and the facts that $\|\omega_x\|_{L^2}\leq  \|\omega_{xx}\|_{L^2}$ and $\|\Omega_x\|_{L^2}\leq  \|\Omega_{xx}\|_{L^2}$. Other terms on the right hand side of (\ref{est-xx1}) can be handled similarly as above. We conclude 
\begin{equation}\notag%\label{est-xx2}
\begin{split}
&\frac12\frac{d}{dt}\left(\|\Omega_{xx}\|_{L^2}^2+\|\omega_{xx}\|_{L^2}^2\right)\\
\leq& \ c_0(1+|a|)\left(\|\Omega_{x}\|_{L^2}+\|\omega_{x}\|_{L^2}\right)\left(\|\Omega_{xx}\|_{L^2}^2+\|\omega_{xx}\|_{L^2}^2\right),
\end{split}
\end{equation}
which immediately gives, by Gr\"onwall's inequality
\begin{equation}\label{est-xx3}
\begin{split}
&\left(\|\Omega_{xx}(t)\|_{L^2}^2+\|\omega_{xx}(t)\|_{L^2}^2\right)\\
\leq& \left(\|\Omega_{xx}(0)\|_{L^2}^2+\|\omega_{xx}(0)\|_{L^2}^2\right) e^{\int_0^t 2c_0(1+|a|)\left(\|\Omega_{x}(\tau)\|_{L^2}+\|\omega_{x}(\tau)\|_{L^2}\right)\,d\tau}.
\end{split}
\end{equation}
Combining (\ref{est-xx3}) with the assumption that $\Omega, \omega\in C([0,T]; H^1)$, it follows that 
\begin{equation}\notag
\sup_{0\leq t\leq T}\left(\|\Omega(t)\|_{H^2}+\|\omega(t)\|_{H^2}\right)<\infty.
\end{equation}

\cbdu

\bigskip

\section{Pure transport case}
\label{sec-transport}

In this section we prove Theorem \ref{thm-tran}. According to Theorem \ref{thm-existence}, there exists a unique solution $(\Omega(t),\omega(t))$ of (\ref{mhd-tran}) on $[0,T]$ for some $T>0$. In view of Theorem \ref{thm-criterion}, in order to show the global existence, it is sufficient to prove 
\[\int_0^T\left(\|H\Omega (t)\|_{L^\infty}+\|H\omega(t)\|_{L^\infty}\right) \, dt<\infty \ \ \ \mbox{for all} \ \ T>0.\]
On the other hand, due to the boundedness of Hilbert transform, we have 
\begin{equation}\notag
\begin{split}
\|H\Omega(t)\|_{L^\infty}\leq c_0\|H\Omega(t)\|_{C^\beta} \leq c_0\|\Omega(t)\|_{C^\beta},\\
\|H\omega(t)\|_{L^\infty}\leq c_0\|H\omega(t)\|_{C^\beta} \leq c_0\|\omega(t)\|_{C^\beta}
\end{split}
\end{equation}
for $\beta\in(0,1)$. As a consequence, we only need to prove:
\begin{Proposition}\label{prop1}
Assume $\Omega_0,\omega_0\in H^1(S^1)$. Let $(\Omega(t),\omega(t))$ be the solution of (\ref{mhd-tran}) with initial data $(\Omega_0,\omega_0)$ on $[0,T]$. Then there exists $\beta_1, \beta_2\in (0,1)$ such that
\begin{equation}\label{est-infty}
\sup_{0\leq t\leq T}\left(\|\omega(t)\|_{C^{\beta_1}}+\|\Omega(t)\|_{C^{\beta_2}}\right)<\infty.
\end{equation}
\end{Proposition}
\pf
Recall the equations satisfied by $(\Omega, \omega)$,
\begin{equation}\notag
\begin{split}
 \Omega_t+ m\Omega_x=&\ 0,\\
 \omega_t+ p\omega_x=&\ 0.
\end{split}
\end{equation}
Consider the characteristics $X_t(x)$ and $Y_t(x)$ satisfying 
\begin{equation}\label{chara1}
\frac{d}{dt} X_t= p(t, X_t(\xi)), \ \ \ X_0(\xi)=\xi, 
\end{equation}
\begin{equation}\label{chara2}
\frac{d}{dt} Y_t= m(t, Y_t(\xi)), \ \ \ Y_0(\xi)=\xi, 
\end{equation}
such that 
\begin{equation}\label{traject}
\Omega(t, Y_t(x))=\Omega_0(x), \ \ \omega(t, X_t(x))=\omega_0(x).
\end{equation}
We notice that there exists a unique solution $X_t(x)$ to the Cauchy problem (\ref{chara1}) and a unique solution $Y_t(x)$ to (\ref{chara2}). Indeed, since $\Omega(t), \omega(t)\in H^1(S^1)\subset C^{\frac12}(S^1)$ and the Hilbert transform is bounded on $C^\beta$, we have 
\begin{equation}\notag
\begin{split}
\|p(t)\|_{C^{1,\frac12}}\leq &\ c_0\|\Omega(t)\|_{C^{1,\frac12}}<\infty,\\
\|m(t)\|_{C^{1,\frac12}}\leq &\ c_0\|\omega(t)\|_{C^{1,\frac12}}<\infty.\\
\end{split}
\end{equation}
Hence, $p$ and $m$ are Lipschitz in time. Thus, the standard ordinary differential equation theory implies existence and uniqueness of solution to (\ref{chara1}) and (\ref{chara2}).

Denote the inverse (backward) trajectory of $X_t(x)$ and $Y_t(x)$ by $q_1(t,x)=X_t^{-1}(x)$ and $q_2(t,x)=Y_t^{-1}(x)$, respectively. Note that $q_1(t,x)$ and $q_2(t,x)$ satisfy respectively,
\begin{equation}\label{chara3}
\partial_t q_1= -p(t, q_1(t,x)), \ \ \ q_1(0,x)=x, 
\end{equation}
\begin{equation}\label{chara4}
\partial_t q_2= -m(t, q_2(t,x)), \ \ \ q_2(0,x)=x. 
\end{equation}

We claim that $p$ and $m$ satisfy the estimate
\begin{equation}\label{est-pm}
\begin{split}
|p(t,x)-p(t,y)|\leq F(|x-y|), \ \ \ x,y\in[-\pi, \pi], \\
|m(t,x)-m(t,y)|\leq G(|x-y|), \ \ \ x,y\in[-\pi, \pi],
\end{split}
\end{equation}
with 
\begin{equation}\label{F}
F(s)=
\begin{cases}
c_0\|\Omega_0\|_{L^\infty}s(1-\log s), \ \ \ 0\leq s\leq 1,\\
c_0\|\Omega_0\|_{L^\infty}, \ \ \ s> 1,
\end{cases}
\end{equation}
and 
\begin{equation}\label{G}
G(s)=
\begin{cases}
c_0\|\omega_0\|_{L^\infty}s(1-\log s), \ \ \ 0\leq s\leq 1,\\
c_0\|\omega_0\|_{L^\infty}, \ \ \ s> 1,
\end{cases}
\end{equation}
for a universal constant $c_0>0$.
We only need to show one of them, for instance, the estimate for $p$. Recall that, by (\ref{H2})
\begin{equation}\notag
p_x(t,x)=H\Omega=\frac{1}{2\pi} P. V. \int_{-\pi}^{\pi}\Omega(t,y)\cot\left(\frac{x-y}{2}\right)\, dy.
\end{equation}
Hence, we have 
\begin{equation}\notag
p(t,x)=\frac{1}{\pi} P. V. \int_{-\pi}^{\pi}\Omega(t,y)\log \left|\sin\left(\frac{x-y}{2}\right)\right|\, dy.
\end{equation}
Without loss of generality, we take $x,y\in(-\pi, \pi)$ such that $-\pi<x<y<\pi$ and 
$\delta= y-x$. We split the interval $[-\pi, \pi]$ into subintervals 
\begin{equation}\notag
I_1=[-\pi, x-\frac{\delta}{2}), \ \ I_2=[x-\frac{\delta}{2}, x+\frac{\delta}2), \ \ I_3=[x+\frac{\delta}2, y+\frac{\delta}2), \ \ 
I_4=[y+\frac{\delta}2, \pi].
\end{equation}
In the case of $x-\frac{\delta}2\leq -\pi$ or $y+\frac{\delta}2> \pi$, we treat $I_1$ or $I_4$ as an empty set. In order to prove the estimate on $p$ in (\ref{est-pm}), we proceed as
\begin{equation}\notag
\begin{split}
|p(t,x)-p(t,y)|=&\left| \frac{1}{\pi} P.V. \int_{-\pi}^{\pi} \Omega(t,z)\left(\log\left|\sin\left(\frac{x-z}{2}\right)\right|-\log\left|\sin\left(\frac{y-z}{2}\right)\right|\right) \,dz \right|\\
\leq& \left| \frac{1}{\pi} P.V. \int_{I_1} \Omega(t,z)\left(\log\left|\sin\left(\frac{x-z}{2}\right)\right|-\log\left|\sin\left(\frac{y-z}{2}\right)\right|\right) \,dz \right|\\
&+ \left| \frac{1}{\pi} P.V. \int_{I_2}\cdot\cdot\cdot \,dz \right|
+ \left| \frac{1}{\pi} P.V. \int_{I_3} \cdot\cdot\cdot\,dz \right|
+ \left| \frac{1}{\pi} P.V. \int_{I_4} \cdot\cdot\cdot\,dz \right|.
\end{split}
\end{equation}
The second term on the right hand side can be estimated as
\begin{equation}\notag
\begin{split}
& \left| \frac{1}{\pi} P.V. \int_{I_2} \Omega(t,z)\left(\log\left|\sin\left(\frac{x-z}{2}\right)\right|-\log\left|\sin\left(\frac{y-z}{2}\right)\right|\right) \,dz \right|\\
\leq& c_0\|\Omega(t)\|_{L^\infty} P.V. \int_{x-\frac{\delta}{2}}^{x-\frac{\delta}{2}} \left|\log|x-z|\right|+\left|\log|y-z|\right| \,dz \\
\leq& c_0\|\Omega(t)\|_{L^\infty} \delta\left(1+|\log \delta|\right)\\
\leq& 
\begin{cases}
c_0\|\Omega_0\|_{L^\infty} \delta\left(1-\log \delta\right), \ \ \ 0<\delta<1\\
c_0\|\Omega_0\|_{L^\infty}, \ \ \ \delta\geq 1.
\end{cases}
\end{split}
\end{equation}
The integrals on $I_1$, $I_3$ and $I_4$ can be estimated similarly. The estimate for $m$ in (\ref{est-pm}) can be established in an analogous way.

In view of (\ref{chara1})-(\ref{chara2}) and (\ref{est-pm}), we have 
\begin{equation}\label{est-qxy1}
\begin{split}
\partial_t\left|q_1(t,x)-q_1(t,y)\right|\leq&\ F(\left|q_1(t,x)-q_1(t,y)\right|), \\
\partial_t\left|q_2(t,x)-q_2(t,y)\right|\leq&\ G(\left|q_2(t,x)-q_2(t,y)\right|). \\
\end{split}
\end{equation}
Denote $\beta_1(t)=e^{-c_0\|\Omega_0\|_{L^\infty} t}$ and $\beta_2(t)=e^{-c_0\|\omega_0\|_{L^\infty} t}$. For fixed $x$ and $y$ with $|x-y|<1$, define
\begin{equation}\notag
z_1(t)=
\begin{cases}
|x-y|^{\beta_1(t)}e^{1-\beta_1(t)}, \ \ \ \ 0\leq t<t_0,\\
1+c_0\|\Omega_0\|_{L^\infty}(t-t_0),  \ \ \ \ t\geq t_0,
\end{cases}
\end{equation}
where $t_0$ is such that $|x-y|^{\beta_1(t_0)}e^{1-\beta_1(t_0)}=1$. Note that $\beta_1(0)=1$ and $z_1(0)=|x-y|<1$. Hence, $z_1(t)$ is well-defined on $[0,\infty)$. One can verify that $z_1(t)$ is the solution of the differential equation 
\begin{equation}\notag
\partial_t z=F(z), \ \ \ z(0)=|x-y|.
\end{equation}
Combining with the first inequality of (\ref{est-qxy1}), we conclude 
\begin{equation}\label{est-qxy2}
|q_1(t,x)-q_1(t,y)|\leq z_1(t). 
\end{equation}
Similarly, we define 
\begin{equation}\notag
z_2(t)=
\begin{cases}
|x-y|^{\beta_2(t)}e^{1-\beta_2(t)}, \ \ \ \ 0\leq t<t_0,\\
1+c_0\|\omega_0\|_{L^\infty}(t-t_0),  \ \ \ \ t\geq t_0,
\end{cases}
\end{equation}
with $t_0$ such that $|x-y|^{\beta_2(t_0)}e^{1-\beta_2(t_0)}=1$. Analogously, using the second inequality of (\ref{est-qxy1}), we infer 
\begin{equation}\label{est-qxy3}
|q_2(t,x)-q_2(t,y)|\leq z_2(t). 
\end{equation}

We are ready to show (\ref{est-infty}). Noticing that $\omega(t, x)=\omega(0, X_t^{-1}(x))$, we deduce
\begin{equation}\notag
\begin{split}
\left|\omega(t, x)-\omega(t, y)\right|=&\ \left|\omega(0, X_t^{-1}(x))-\omega(0, X_t^{-1}(y))\right|\\
=&\left|\int_{X_t^{-1}(y)}^{X_t^{-1}(x)} \omega_{0,x}(\zeta)\, d\zeta\right|\\
\leq&\ c_0\|\omega_{0,x}\|_{L^2}\left|X_t^{-1}(x)-X_t^{-1}(y)\right|^{\frac12}\\
\leq&\ c_0\|\omega_{0,x}\|_{L^2}\left|q_1(t,x)-q_1(t,y)\right|^{\frac12}\\
\end{split}
\end{equation}
where mean value theorem and H\"older's inequality were applied. As a consequence, we conclude 
\begin{equation}\notag
\sup_{0\leq t\leq T}\|\omega(t)\|_{C^{\beta_1}}<\infty.
\end{equation}
thanks to (\ref{est-qxy2}). Analogously, we can show
\begin{equation}\notag
\sup_{0\leq t\leq T}\|\Omega(t)\|_{C^{\beta_2}}<\infty.
\end{equation}
It completes the proof of the proposition.

\cbdu

\bigskip

\section{Numerical simulations}
\label{sec-numerics}

In this section, we perform some numerical study for the 1D model (\ref{mhd-1d}) of MHD. For convenience, we recall (\ref{mhd-1d}) here,
%We consider periodic solutions of the system
\begin{equation}\label{mhd-1d-6}
\begin{split}
\Omega_t+am \Omega_x-\omega p_x=&\ 0, \\
\omega_t+ap\omega_x - \Omega m_x=&\ 0, \\
p_x=H \Omega, \ \ 
m_x=&\ H \omega,
\end{split} \quad x\in [-\pi, \pi]
\end{equation}
%for $a=1$ and $a=-1$. 
and the Hilbert transform for a periodic function
\begin{equation}\notag
H f=\frac{1}{2\pi} P.V. \int_{-\pi}^{\pi}f(y)\cot\left( \frac{x-y}{2}\right)\, dy.
\end{equation}
%Notice that $p$ and $m$ defined by antiderivatives are not unique. 
As mentioned earlier, in order for $p$ and $m$ to be uniquely defined, 
we can choose the gauge and set them to have either zero mean over the interval $[-\pi, \pi]$ or zero point value at a fixed point, e.g., $p(x_0, t)=m(x_0, t)=0$ for some $x_0$, see \cite{Jia2019}.

We use a Fourier-collocation spectral method for the spatial approximation and a five stage fourth order low storage Runge-Kutta method for time discretization. 
An exponential type filter is used  for stabilization of the spectral method, see \cite{Hesthaven}. 
For a periodic function $f(x)$,  its Hilbert transform can be approximated in spectral method via the following formula: 
\[\widehat{Hf}(k)=-i \mbox{sgn} (k)\hat f(k),\]
where $\hat f(k)$ are  coefficients in  Fourier series of $f(x)$, 
see \cite{Hou, OSW}.
Similarly, for periodic functions $p(x)$ and $\Omega(x)$, the equation $p_x=H \Omega$ can be  approximated in spectral method through the relation
\[ i k \hat p(k)=-i \mbox{sgn} (k)\hat \Omega (k).\]

\medskip

\subsection{Numerical results for the 1D model of MHD}

%Applying the following simple properties of the Hilbert transform 
%\begin{equation}\notag
%\begin{split}
%H(cf)=&\ cHf, \ \ \mbox{for a constant} \ \ c,\\
%H\sin(kx)=&-\cos(kx),\ \ \ H\cos(kx)=\sin(kx),
%\end{split}
%\end{equation}
%we have 
%\begin{equation}\notag
%H\sin(kx+\theta)=-\cos(kx+\theta), \ \ H\cos(kx+\theta)=\sin(kx+\theta).
%\end{equation}
One can check that, for arbitrary constants $A_1$, $A_2$, $\theta_1$, $\theta_2$ and $k$ 
\begin{equation}\notag
\Omega(x)=A_1\sin(kx+\theta_1), \ \ \ \omega(x)=A_2\sin(kx+\theta_2)
\end{equation}
and 
\begin{equation}\notag
\Omega(x)=A_1\cos(kx+\theta_1), \ \ \ \omega(x)=A_2\cos(kx+\theta_2)
\end{equation}
are steady states of system (\ref{mhd-1d-6}). Thus, we choose to consider 
the initial condition
\begin{equation}\label{initial-num}
\Omega_0=\sin(x)+\cos(4x)+5, \ \ \ \omega_0=\sin(2x)+2,
\end{equation}
such that there are two non-zero modes in the Fourier representation of $\Omega_0$.

We conduct simulations for (\ref{mhd-1d-6}) with initial data (\ref{initial-num}) under the following two settings: (i) $a=1$ and (ii) $a=-1$.
% \begin{enumerate}
%\item[(i)] $a=1$; % and $\int_{-\pi }^{\pi} p(t, x)\, dx=\int_{-\pi }^{\pi} m(t, x)\, dx=0$. 
%\item[(ii)] $a=-1$.% and $\int_{-\pi }^{\pi} p(t, x)\, dx=\int_{-\pi }^{\pi} m(t, x)\, dx=0$. 
% \end{enumerate}
In the computation, we take $N=12800$ points in the Fourier-collocation spectral method. 
 
Figure \ref{Fig-case2} shows the numerical results for case (i). The time evolution of $\Omega(t,x)$ and $\omega(t,x)$ are plotted in Figure \ref{Fig-case2}(a) and Figure \ref{Fig-case2}(b), respectively. One can see that $\Omega(t,x)$ and $\omega(t,x)$ are rather smooth. The first order derivative $\Omega_x$ shown in Figure \ref{Fig-case2}(c) seems smooth as well, while $\omega_x$ illustrated in Figure \ref{Fig-case2} (d) develops some mild spines at time $t=4$. However, we observe spines for the second derivatives $\Omega_{xx}$ and $\omega_{xx}$ at larger time in Figure \ref{Fig-case2}(e) and (f).  In particular, there is a notable spine near $x=0$ at $t=4$. Notice that 
\[u=\frac12(p+m), \ \ \ B=\frac12(p-m),\]
and hence
\begin{equation}\label{ub-pm}
\begin{split}
u_x=\frac12(p_x+m_x)=& \frac12 H(\Omega+\omega), B_x=\frac12(p_x-m_x)= \frac12 H(\Omega-\omega),\\
H\Omega=&\ u_x+B_x, \ \ H\omega=u_x-B_x.
\end{split}
\end{equation}
Figure \ref{Fig-case2}(g) shows the time evolution of $\|H\Omega\|_{L^\infty}+\|H\omega\|_{L^\infty}$, while Figure \ref{Fig-case2}(h) shows $\|u_x(t)\|_{L^\infty}$ and $\|B_x(t)\|_{L^\infty}$. We observe oscillations in these graphs and the amplitudes grow slowly in a linear manner. Combined with the regularity criterion (\ref{criterion}), it seems that the solution starting with data (\ref{initial-num}) may not develop singularities at finite time.

%The results for case (ii) are shown in Figure \ref{Fig-case2}. We note that the behaviours of solutions shown in Figure \ref{Fig-case2} and Figure \ref{Fig-case1} are over all similar. The difference is that spine appears for the first order derivative $\omega_x$ at time $t=4$ in Figure \ref{Fig-case2}(d) and the second derivative $\omega_{xx}$ grows more rapidly in Figure \ref{Fig-case2}(f) than in Figure \ref{Fig-case1}(f). This indicates the choice of gauge with either zero mean or fixing one point does not essentially affect the behaviours of solutions. 

For case (ii) with $a=-1$ and $p(t, 0)=m(t, 0)=0$, the results are illustrated in Figure \ref{Fig-case3}. One can see from
Figure \ref{Fig-case3}(a) and (b) that the solution is less smooth compared to the solution in case (i) shown in Figure \ref{Fig-case2}(a) and (b). This suggests that the convection term with a negative sign causes the solution to behave more singularly. Nevertheless, \ref{Fig-case3}(c) and (d) show that the amplitudes of $\|H\Omega\|_{L^\infty}+\|H\omega\|_{L^\infty}$, $\|u_x(t)\|_{L^\infty}$ and $\|B_x(t)\|_{L^\infty}$ grow faster than that of case (i), but remain in a linear growth. Thus one may speculate that solutions of system (\ref{mhd-1d-6}) with $a=-1$ starting from smooth initial data do not develop singularities in finite time.

 \begin{figure}[!htb]
  \subfigure[]{\includegraphics[scale=0.12]{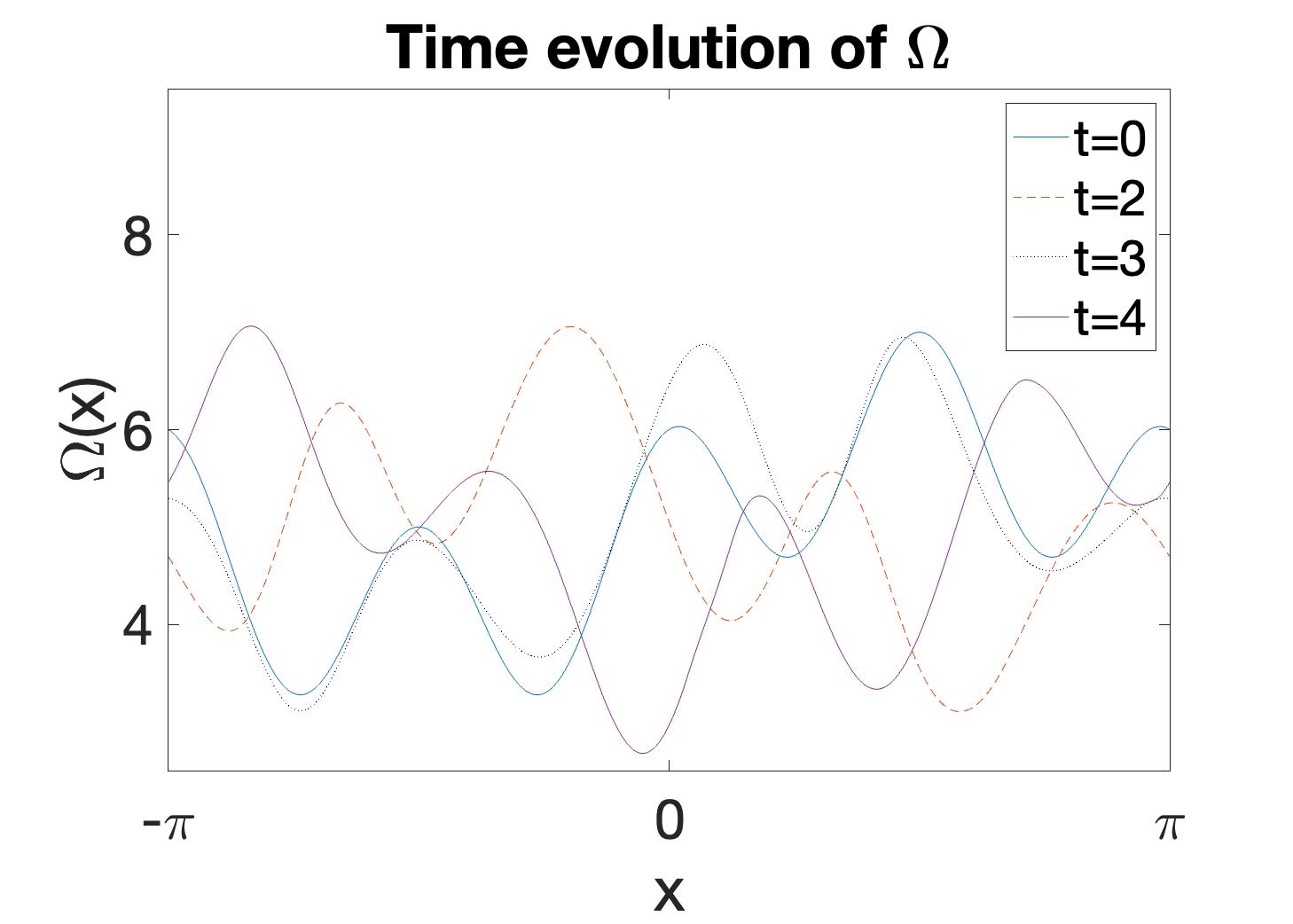} }
  \hspace{-0.2in}
    \subfigure[]{\includegraphics[scale=0.12]{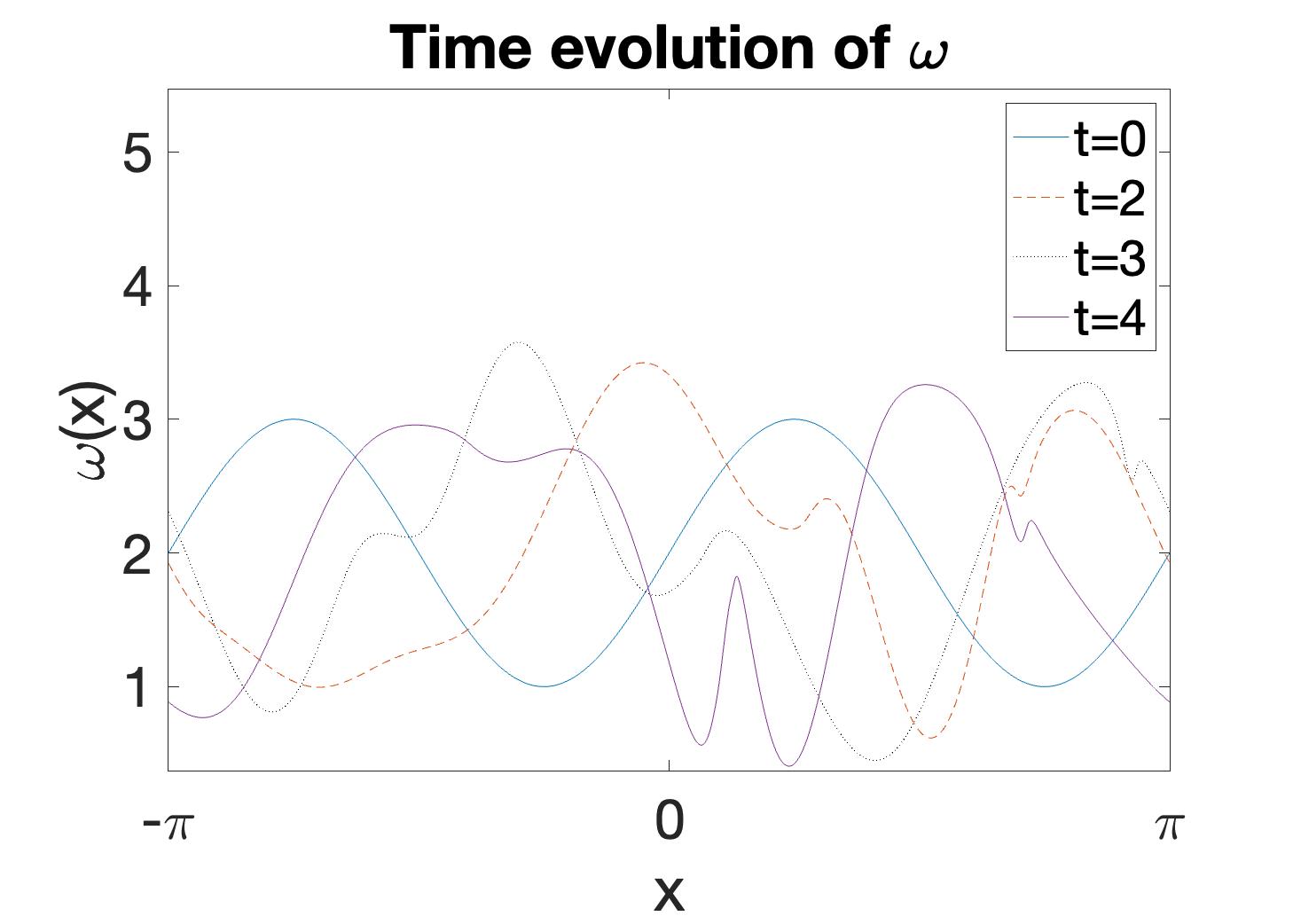} }\\
  \subfigure[]{\includegraphics[scale=0.12]{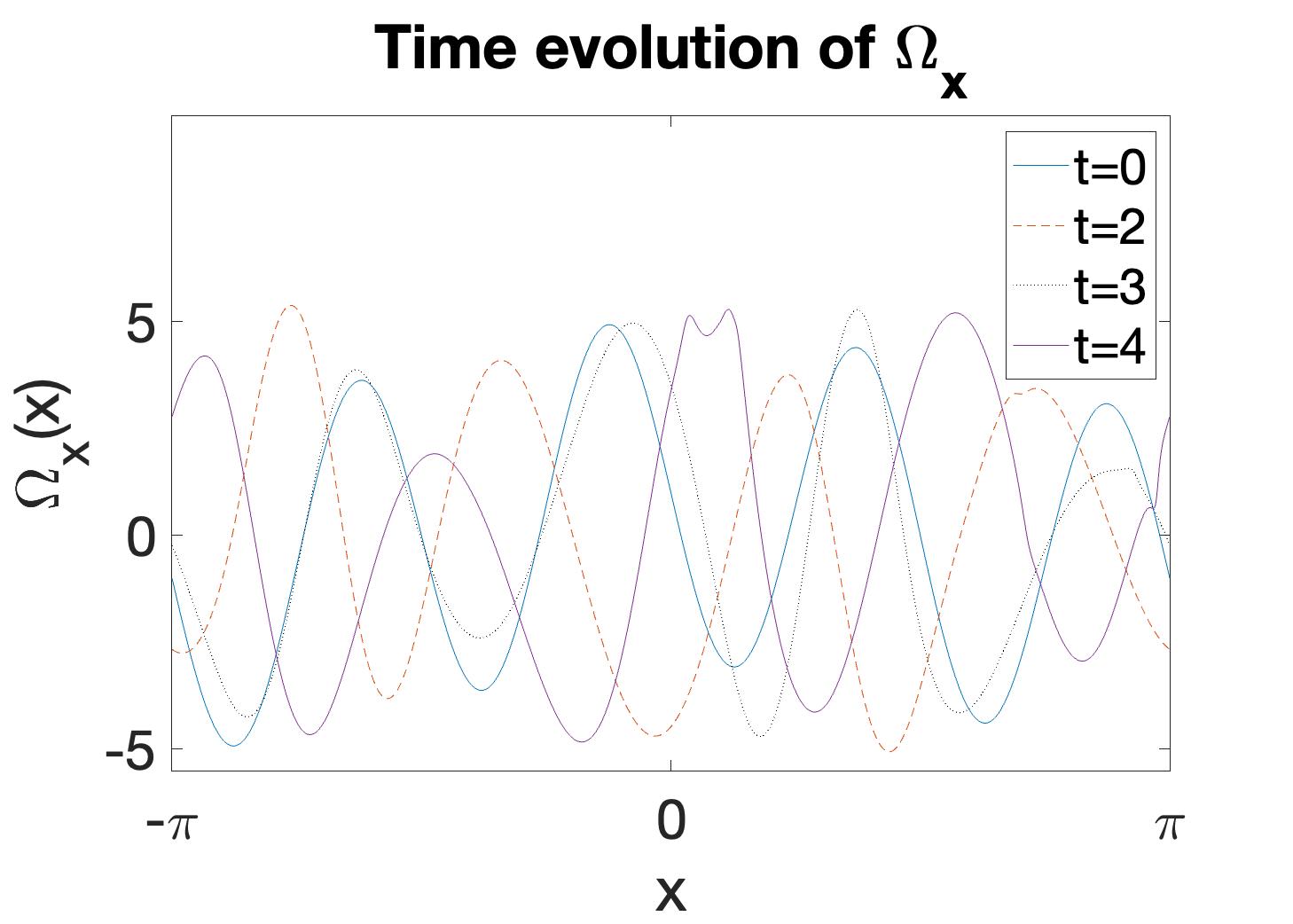} }  \hspace{-0.2in}
    \subfigure[]{\includegraphics[scale=0.12]{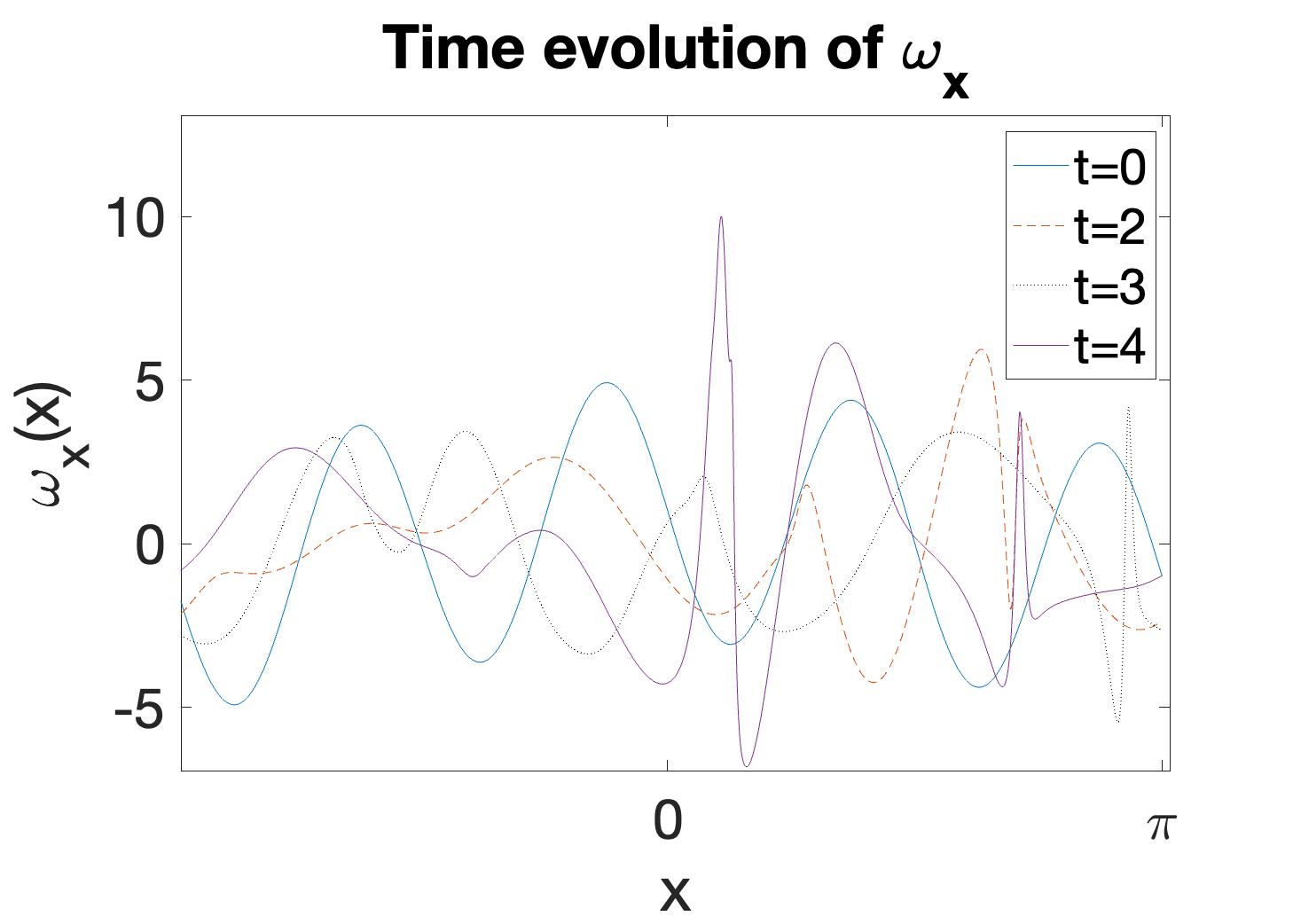} } 
  \\
  \subfigure[]{\includegraphics[scale=0.12]{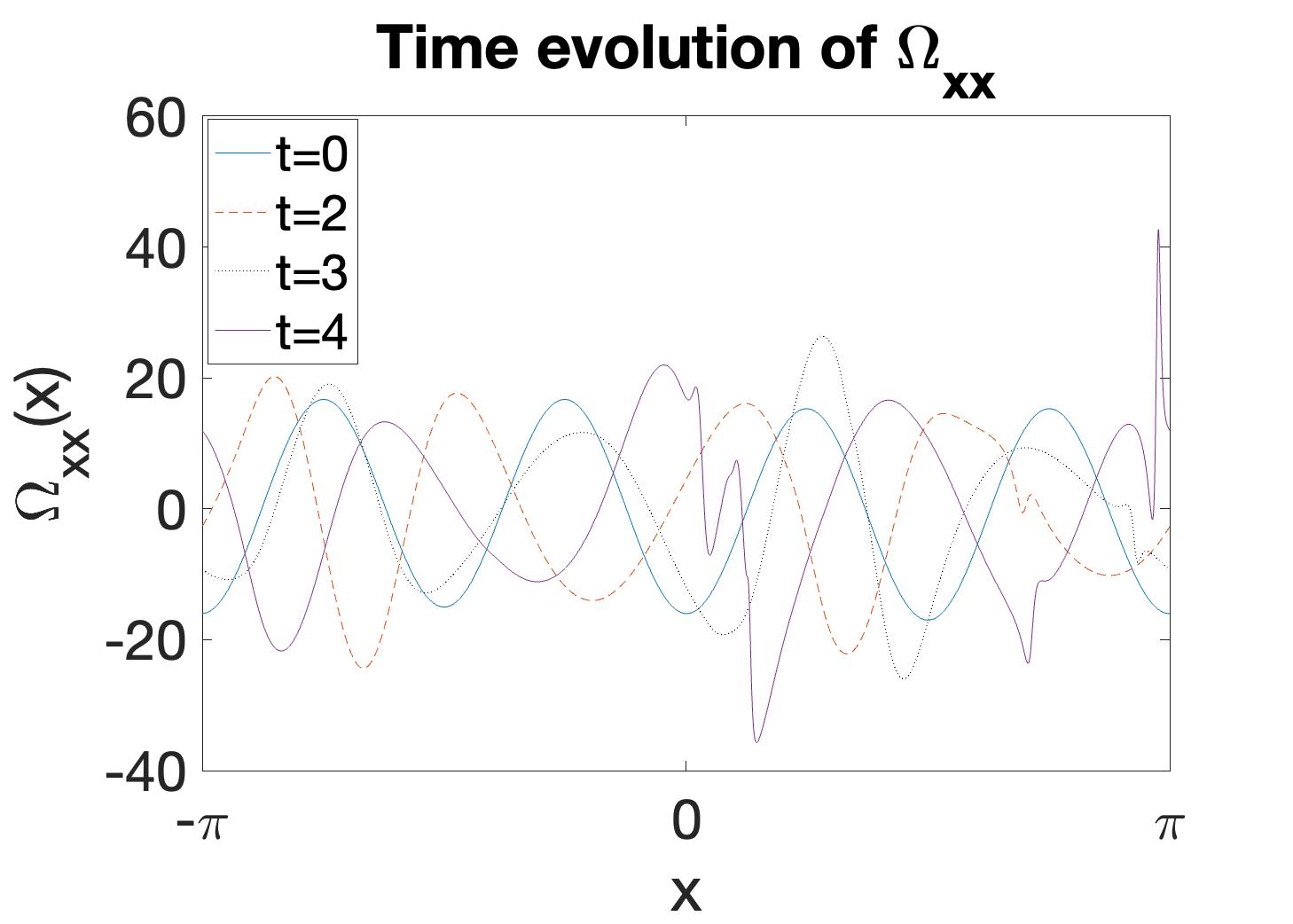} } \hspace{-0.2in}
    \subfigure[]{\includegraphics[scale=0.12]{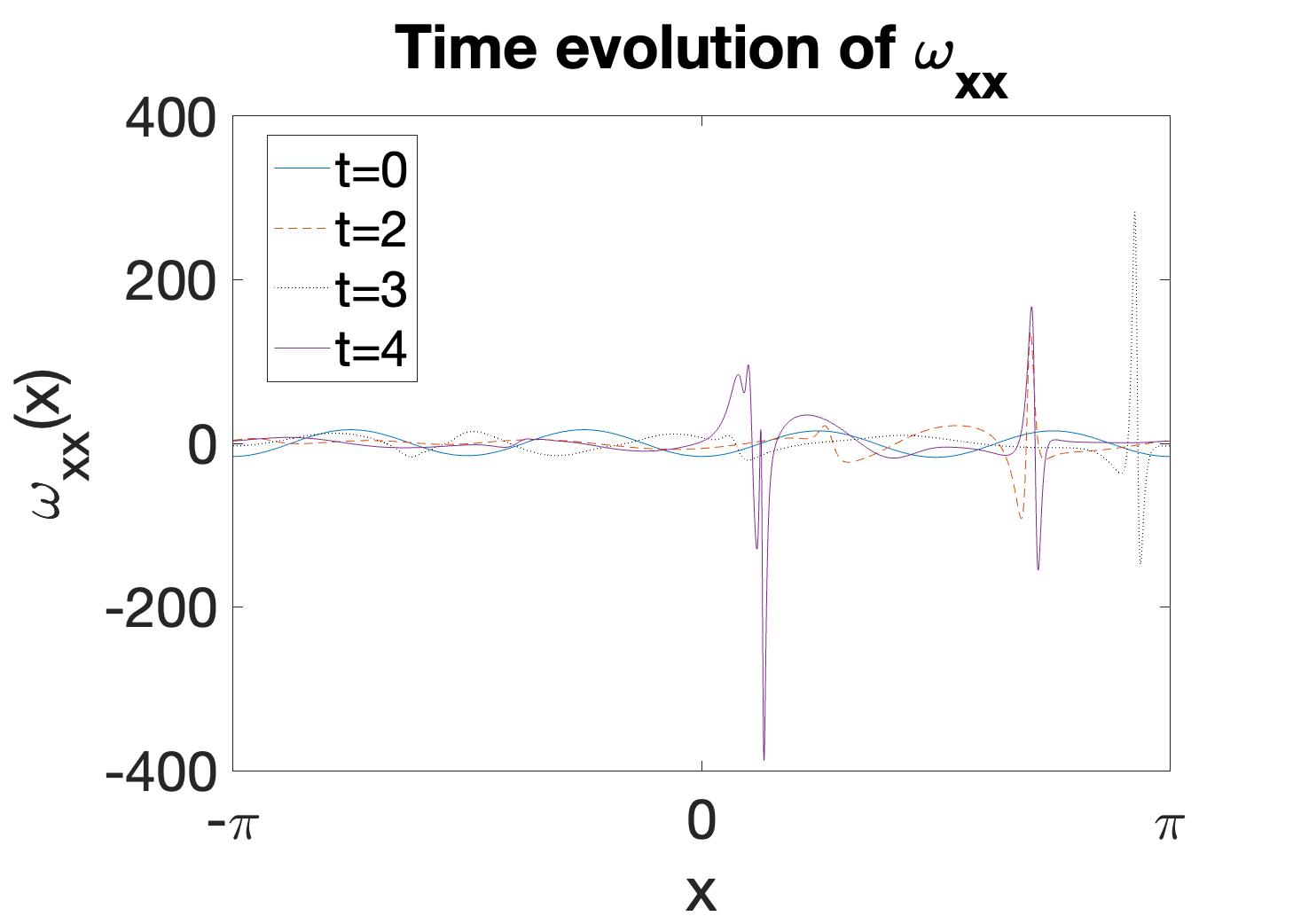} }\\
    \subfigure[]{\includegraphics[scale=0.12]{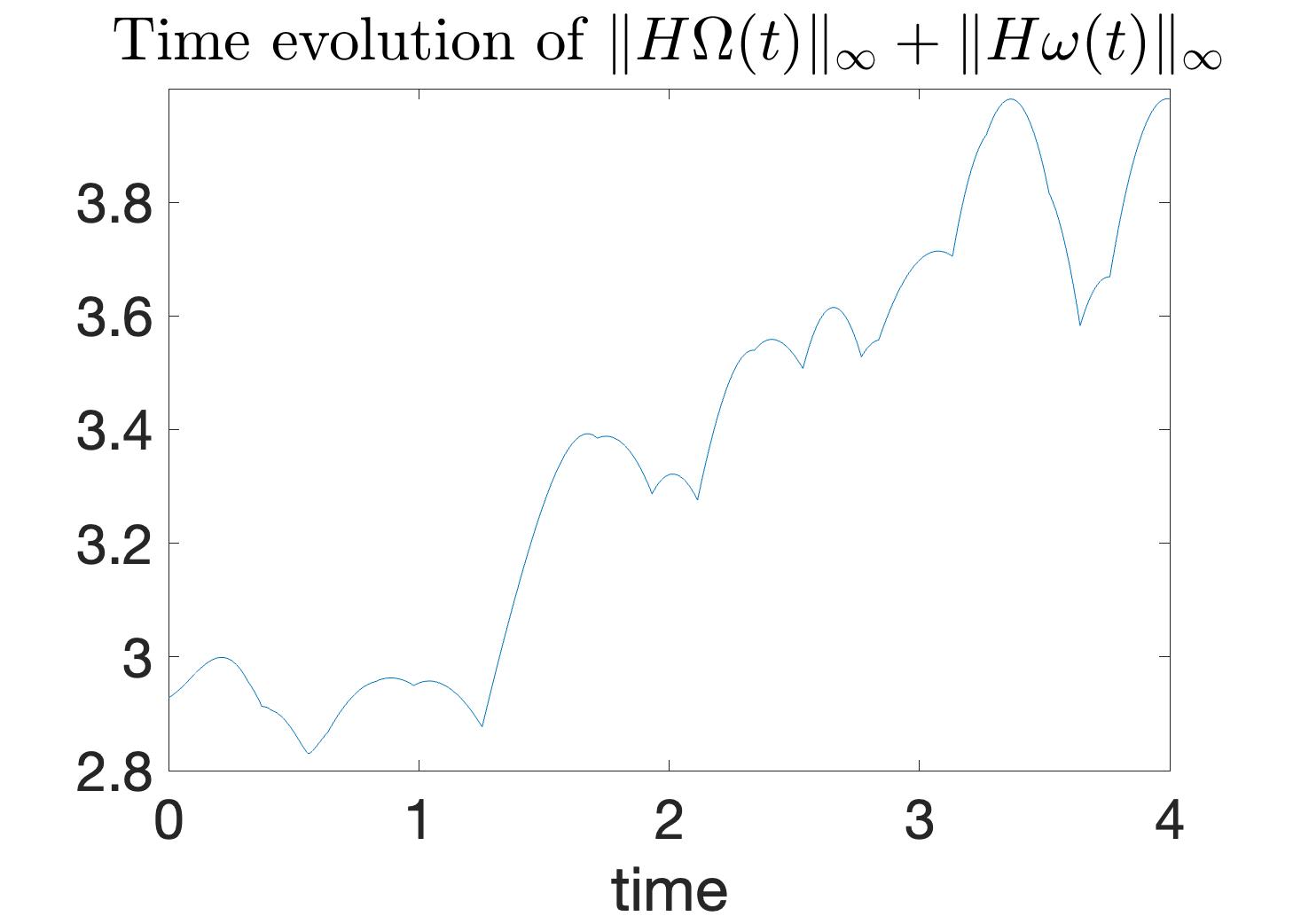} } \hspace{-0.2in}
        \subfigure[]{\includegraphics[scale=0.12]{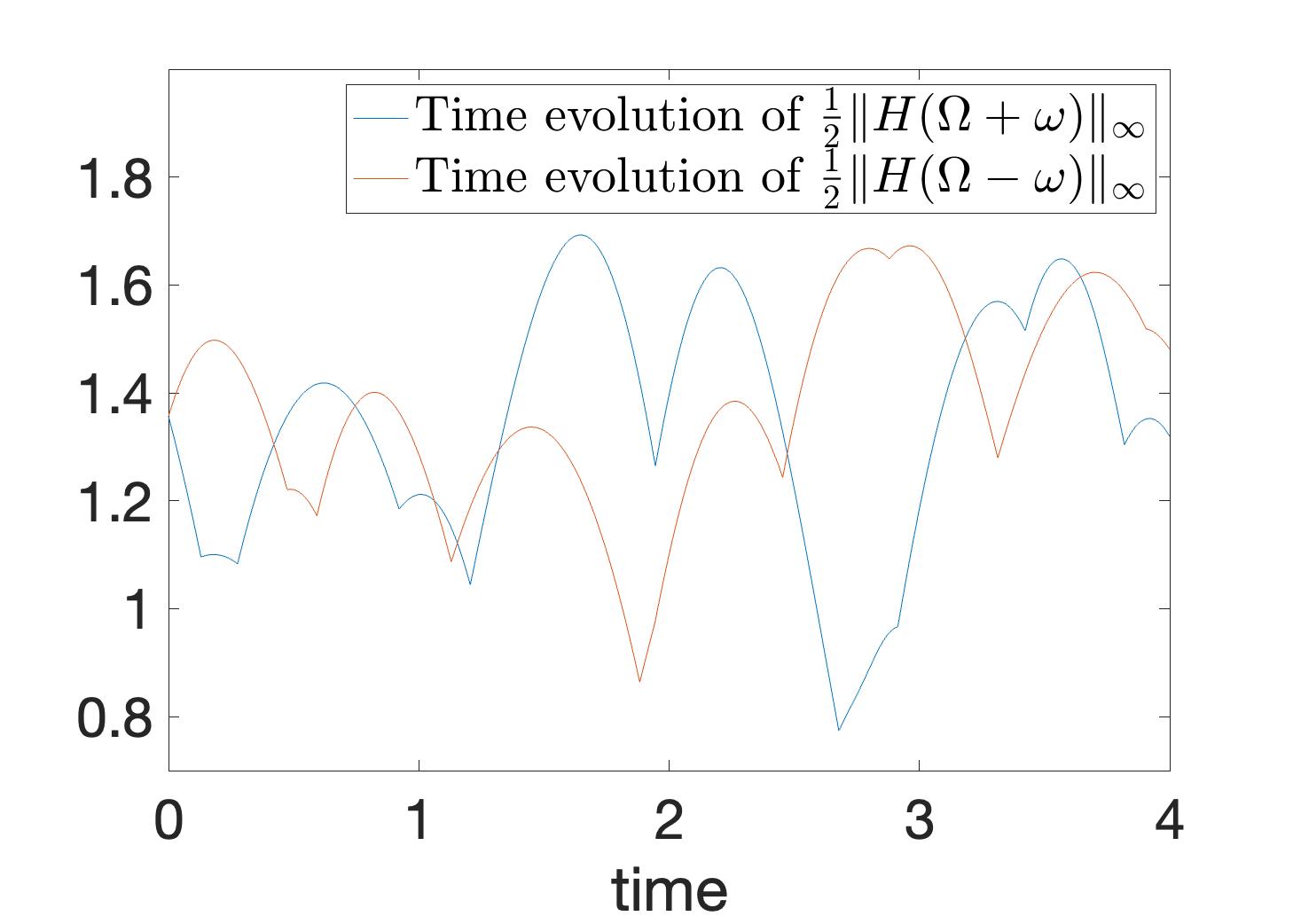} }
  \caption{$a=1$} %{and $\int_{-\pi }^{\pi} p(x, t)\, dx=\int_{-\pi }^{\pi} m(x, t)\, dx=0$.  }
\label{Fig-case2}
 \end{figure}

  \begin{figure}[!htb]
  \subfigure[]{\includegraphics[scale=0.12]{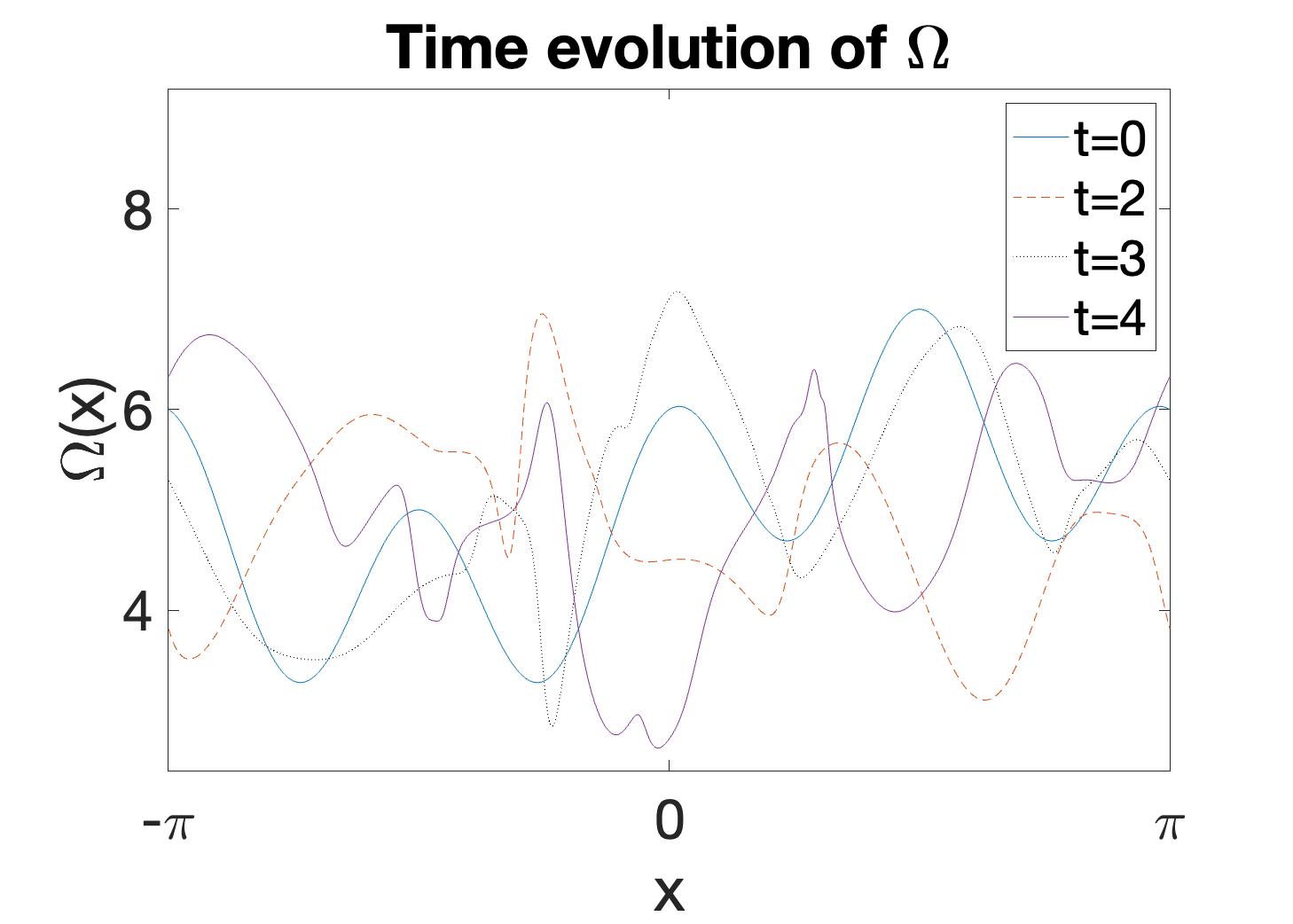} }
  \hspace{-0.2in}
    \subfigure[]{\includegraphics[scale=0.12]{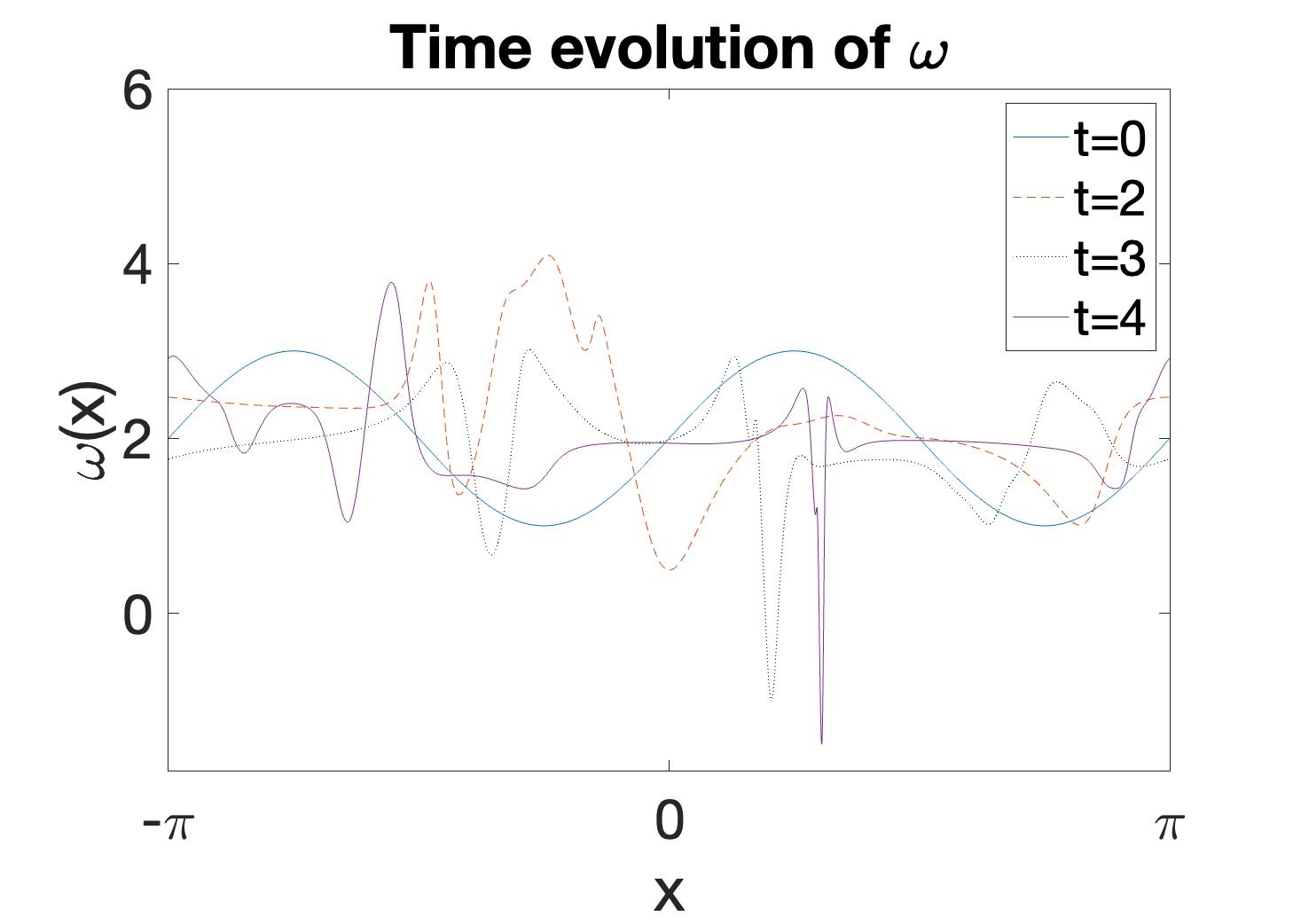} }\\
    \subfigure[]{\includegraphics[scale=0.12]{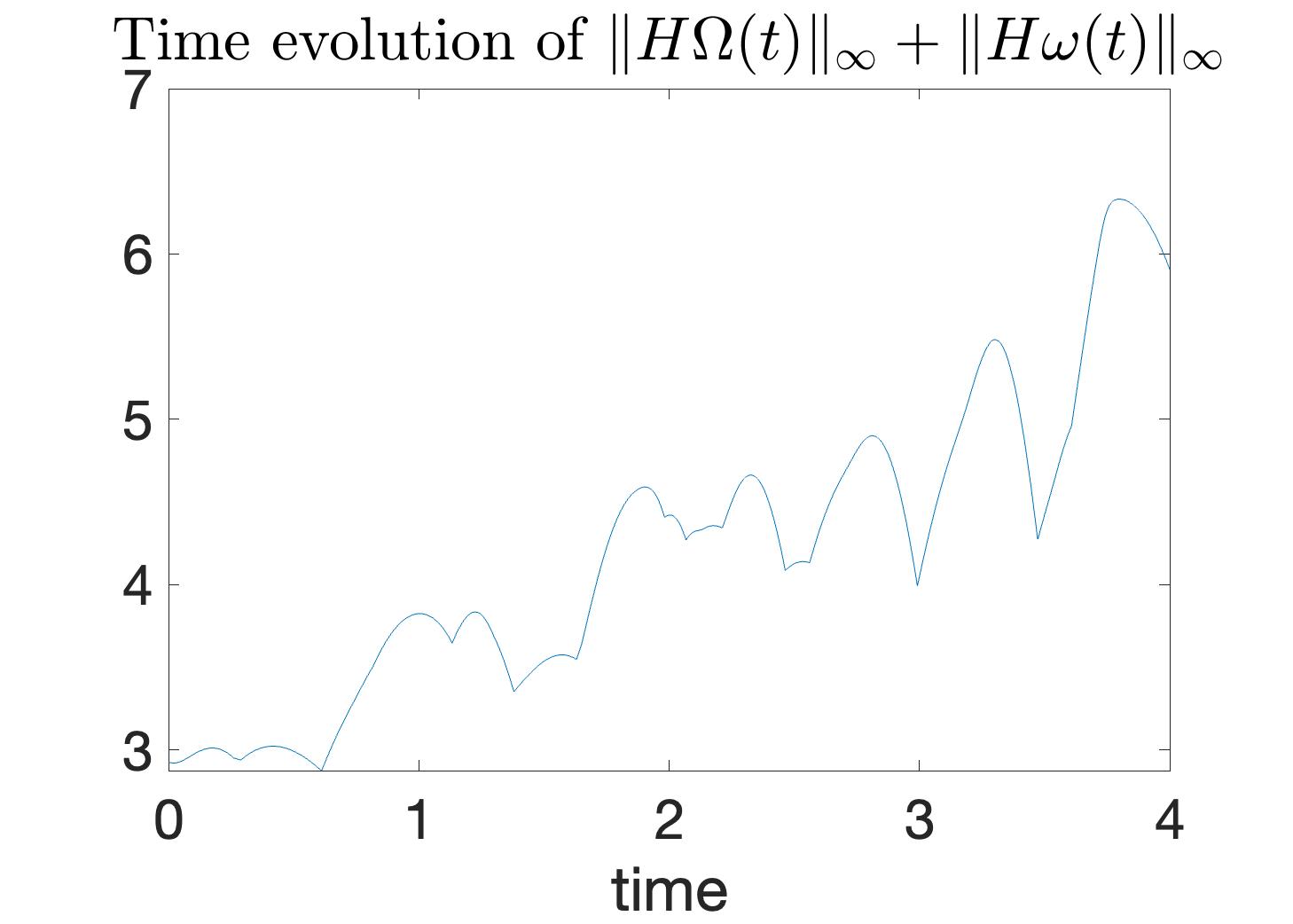} }
      \hspace{-0.2in}
        \subfigure[]{\includegraphics[scale=0.12]{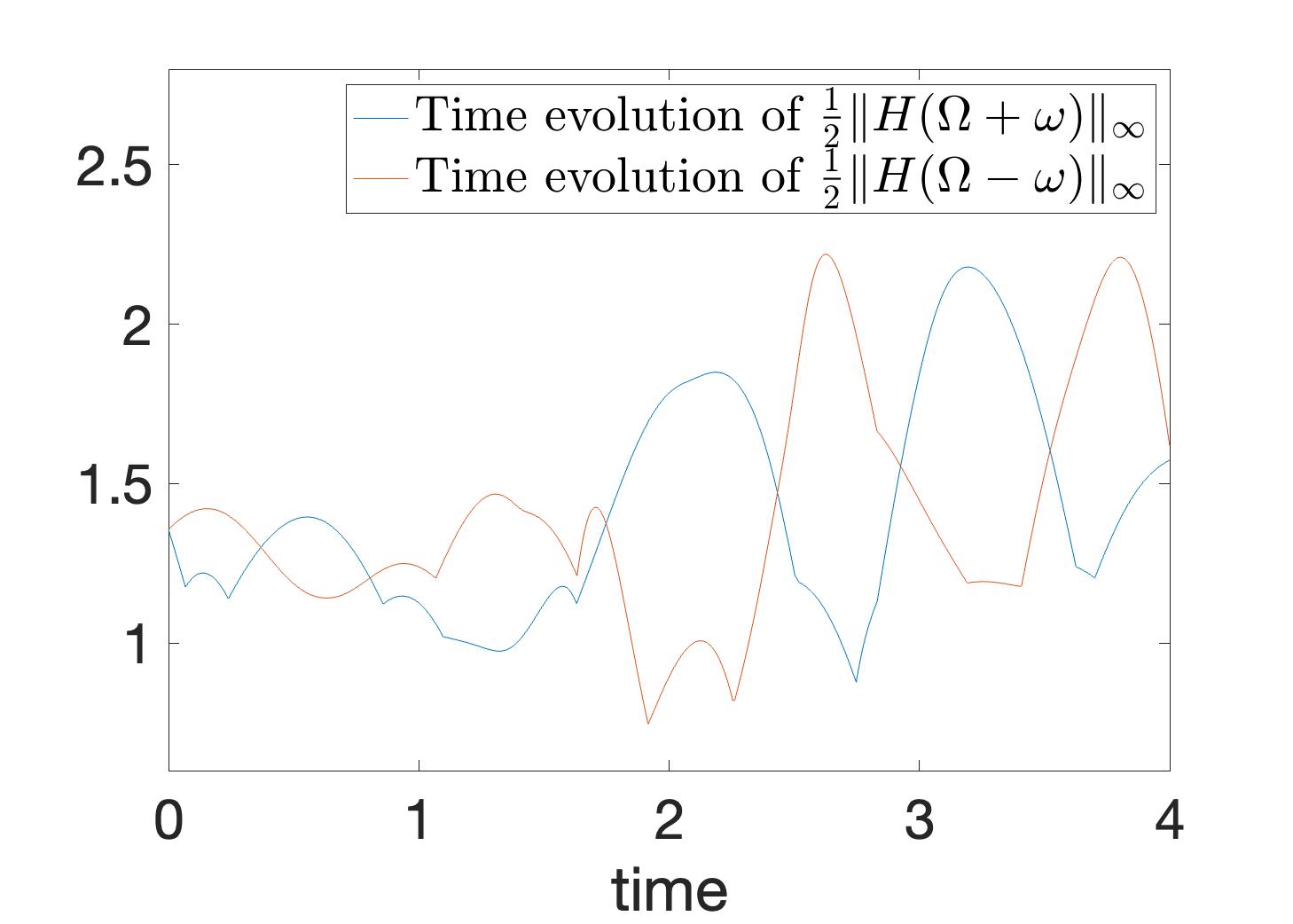} }
  \caption{$a=-1$} %{and $\int_{-\pi }^{\pi} p(x, t)\, dx=\int_{-\pi }^{\pi} m(x, t)\, dx=0$.  }
\label{Fig-case3}
 \end{figure}
 
%\subsection{Numerical results for the De Gregorio model revisited}
 \begin{figure}[!htbp]
  \subfigure[]{\includegraphics[scale=0.12]{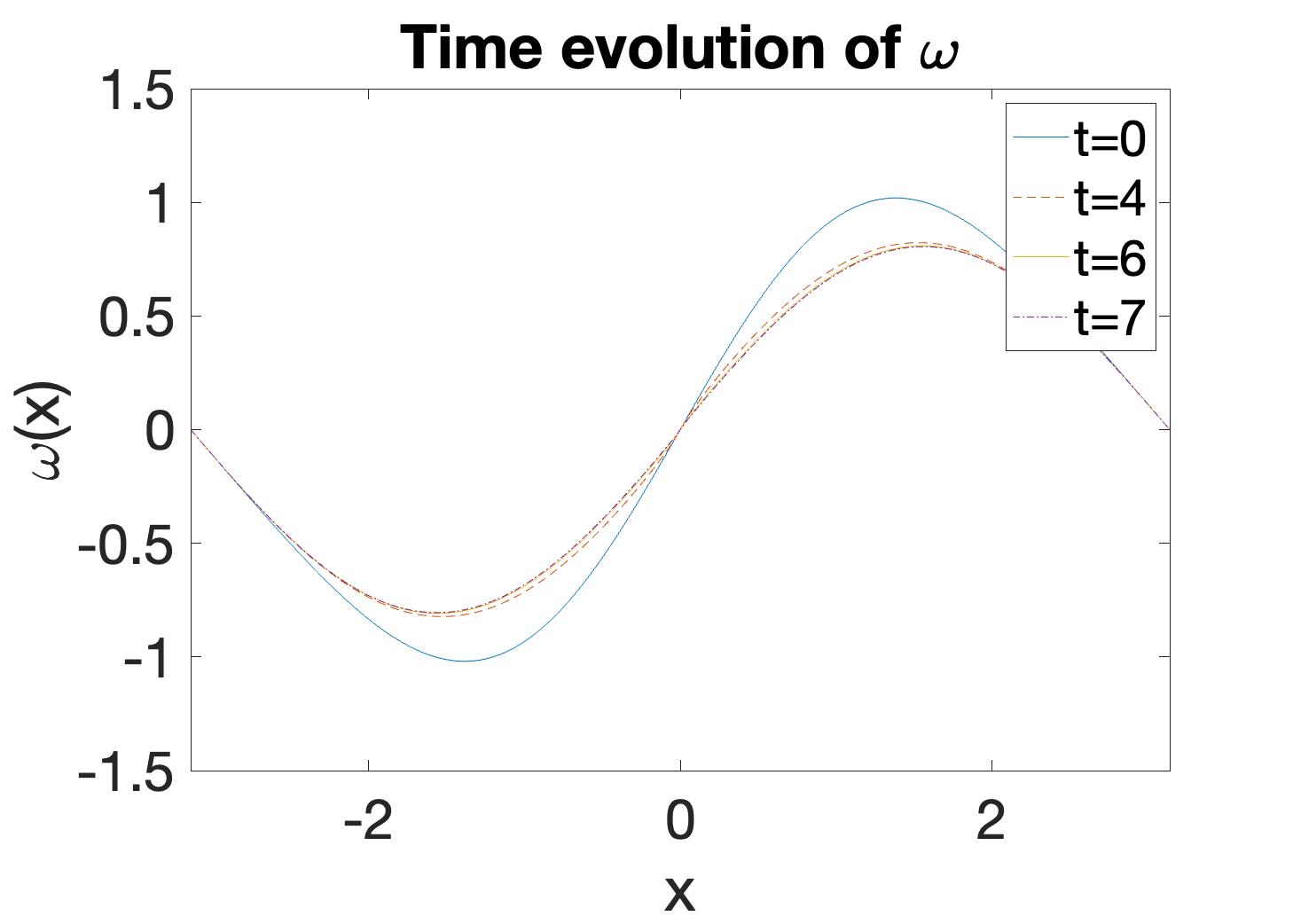} }
  \hspace{-0.2in}
    \subfigure[]{\includegraphics[scale=0.12]{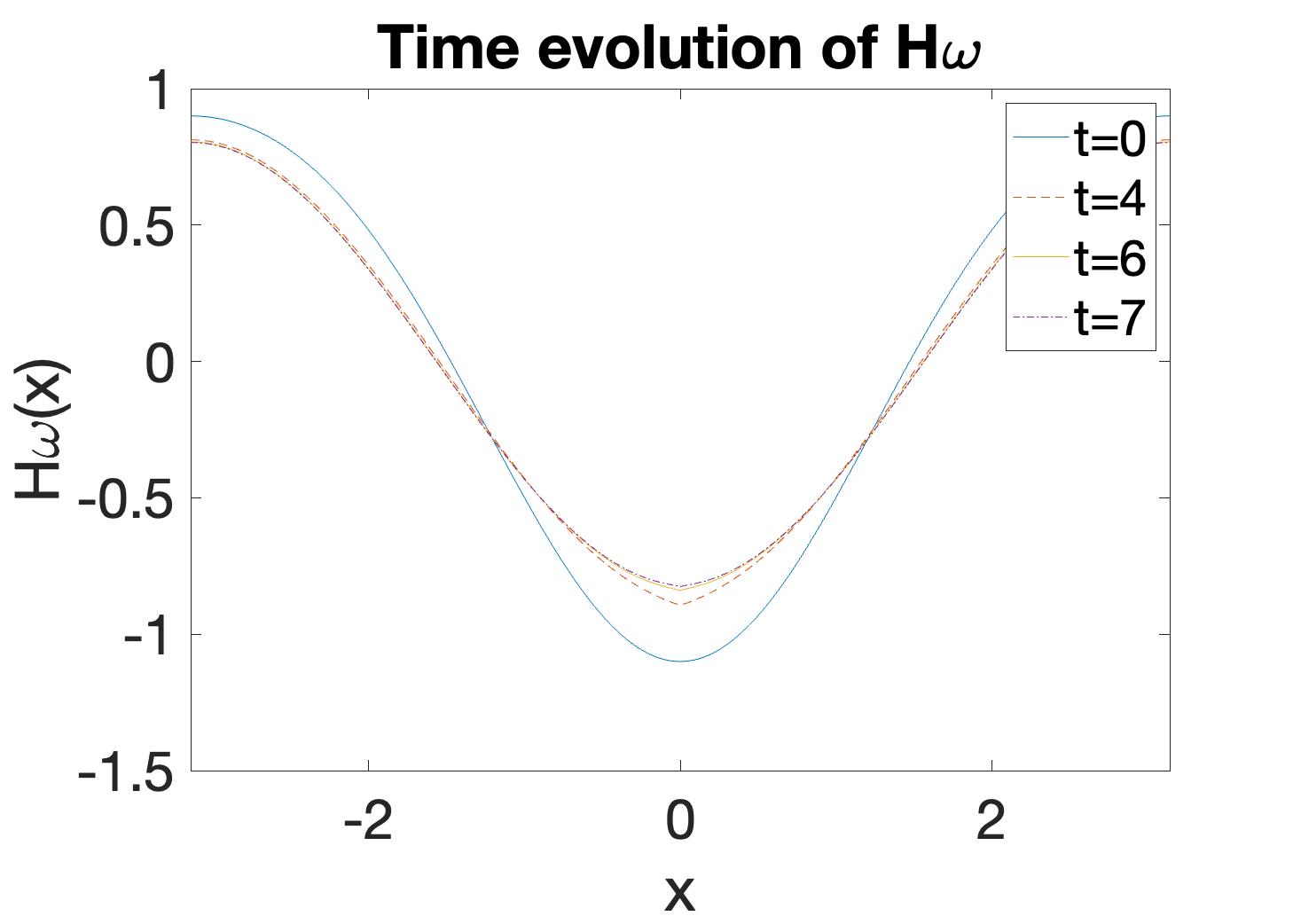} }\\
  \subfigure[]{\includegraphics[scale=0.12]{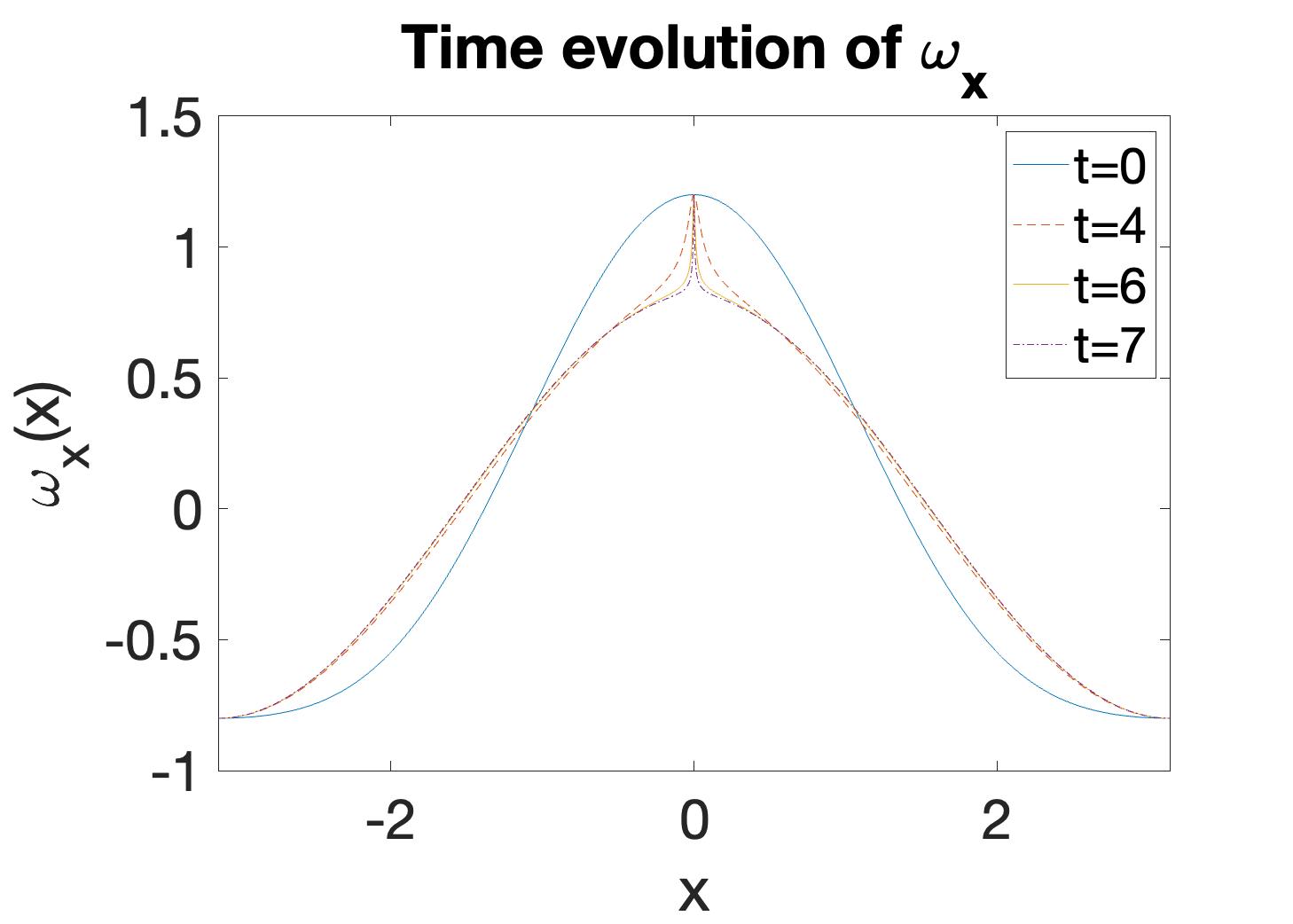} }  \hspace{-0.2in}
    \subfigure[]{\includegraphics[scale=0.12]{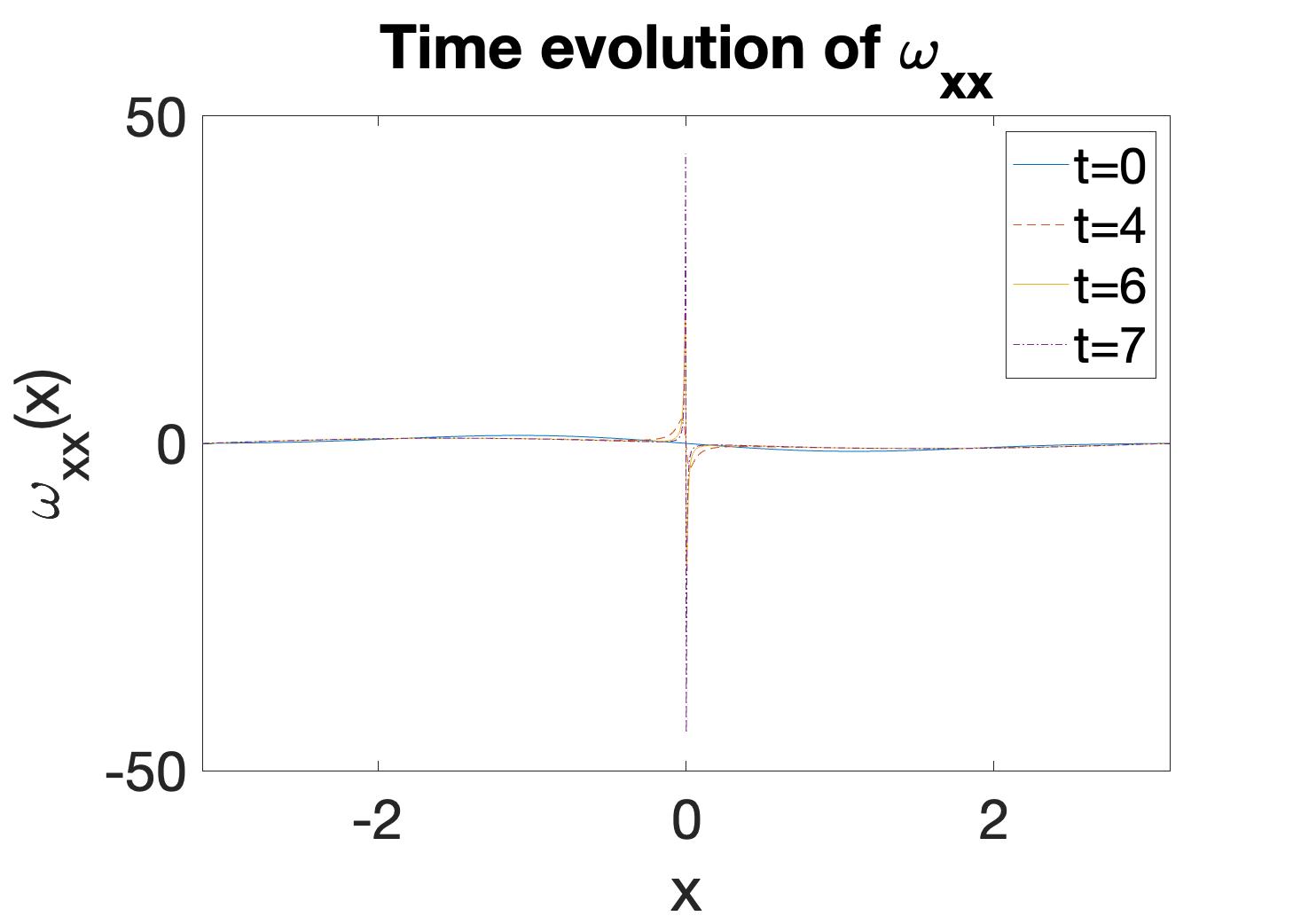} } 
  \\
  \subfigure[]{\includegraphics[scale=0.12]{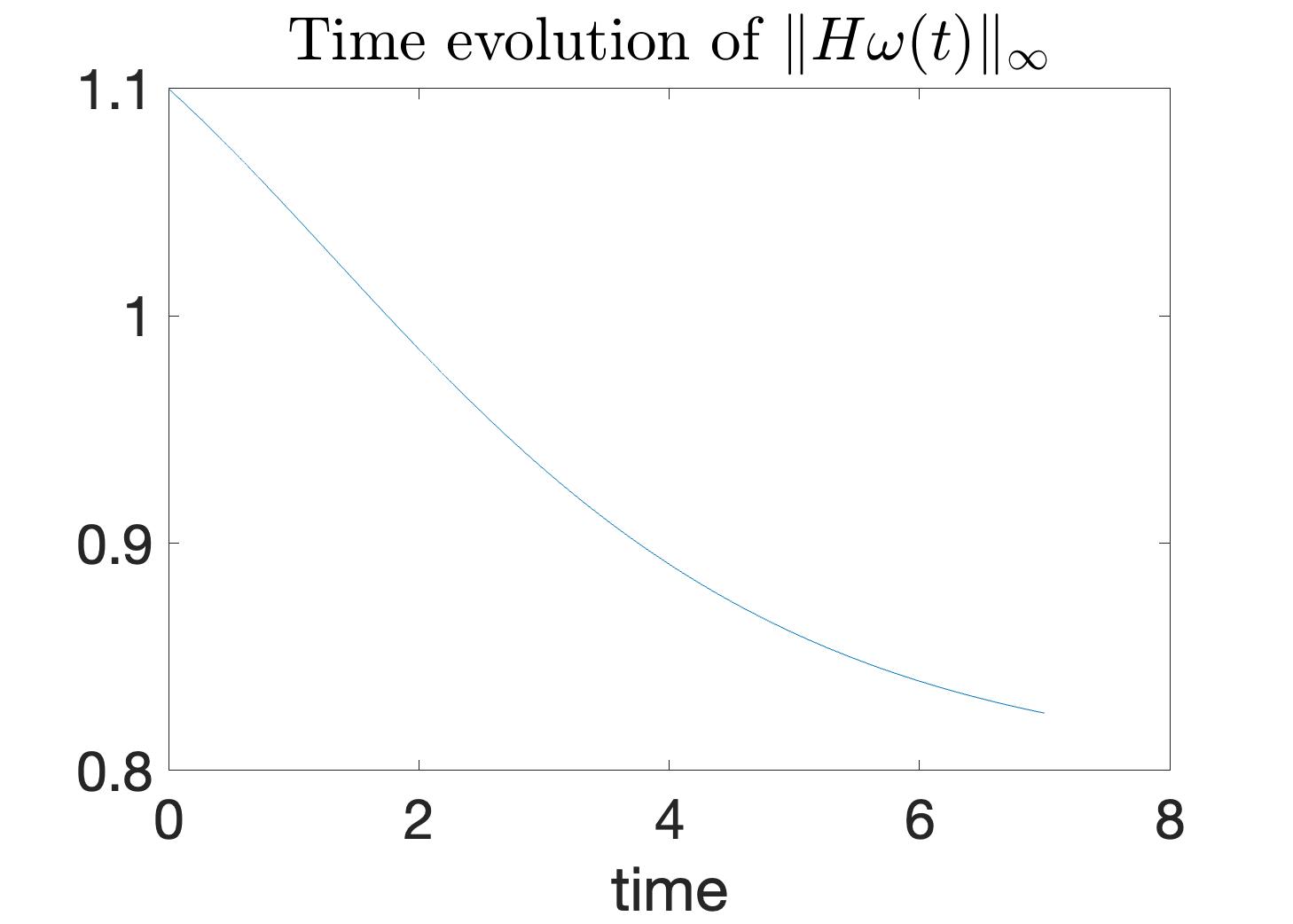} } \hspace{-0.2in}
  \caption{The De Gregorio model.  }
\label{Fig-case4}
 \end{figure}

%Elgindi-Jeong: Numerical simulations done by De Gregorio [9], Okamoto, Sakajo, and Wunsch [20], and others indicate, in fact, that singularity formation for De Gregorio’s model is impossible from smooth data.

\medskip

\subsection{Numerical results for the De Gregorio model revisited}
Numerical simulations for the De Gregorio model (\ref{eq1-deg})-(\ref{eq2-deg}) have been performed in \cite{DeG1, DeG2, OSW} among others. The main information is that singularity formation for this model with smooth initial data is unlikely to happen.

We apply our numerical scheme to (\ref{eq1-deg})-(\ref{eq2-deg}) with the initial data 
\begin{equation}\notag
\omega_0(x)=\sin x+0.1\sin(2x)
\end{equation}
by taking $N=12800$ points in the Fourier-collocation spectral method. The obtained simulations are shown in Figure \ref{Fig-case4}, which recover the numerical results done by Okamoto, Sakajo, and Wunsch \cite{OSW}.
%Plot: $\omega(t,x)$, $H\omega(t,x)$, $\omega_x(t,x)$, $\omega_{xx}(t,x)$ and $\|H\omega\|_{L^\infty}$.
%We note our result is consistent with that done in \cite{OSW}.

We note that $u_x=H\omega$ for the De Gregorio model (\ref{eq1-deg})-(\ref{eq2-deg}) and $u_x=\frac12 H(\Omega+\omega)$ for our 1D MHD model (\ref{mhd-1d-6}), see (\ref{ub-pm}). Comparing Figure \ref{Fig-case2}(h) and Figure \ref{Fig-case4}(e), we observe oscillations of $\|u_x\|_{L^\infty}$ for the 1D MHD model and absence of such oscillations for the pure fluid model. It is reasonable to infer that the interactions between fluid velocity and magnetic field cause such oscillations and more complicated dynamics.

%Comparison of the numerical results for the De Gregorio model and the model for MHD: the interactions between fluid velocity and magnetic field cause oscillations .... 

\bigskip

%\section*{Acknowledgement}
%The authors are grateful for the support of ...

%\Endrefs
\end{document}